  \documentclass[12pt,onecolumn ]{IEEEtran}

\IEEEoverridecommandlockouts

\usepackage{url}
\usepackage{amsmath,amssymb}
\usepackage{amsthm,thmtools}
\usepackage{amsfonts} 

\usepackage[dvips]{epsfig}      
\usepackage{graphicx}
\usepackage{enumerate}
\usepackage{color}                  

\usepackage{verbatim}

\usepackage{amssymb}
\usepackage{amsbsy}
\usepackage{amsmath}
\usepackage{setspace}
\usepackage{url}

\usepackage{subcaption}
\usepackage{float}

 \newcommand{\diag}{\operatorname{diag}}
 
 \newcommand{\tr}{\operatorname{tr}}

\declaretheorem[name={Example},qed={\lower-0.3ex\hbox{$\square$}} ] {Example}

\declaretheorem[name={Definition}  ] {Definition}
\declaretheorem[name={Theorem}  ] {Theorem}
\declaretheorem[name={Lemma}  ] {Lemma}
\declaretheorem[name={Remark}  ] {Remark}
\declaretheorem[name={Corollary}  ] {Corollary}

\declaretheorem[name={Proposition}  ] {Proposition}

\newcommand {\R}{\mathbb R}

\newcommand {\M}{\mathbb M}
\newcommand {\Q}{\mathbb Q}

\newcommand{\be}{\begin{equation}}
\newcommand{\ee}{\end{equation}}

\newcommand{\myc}[1]{{\bf{#1}}}  

  \onehalfspace

\newcommand{\model}{CVDDS}

\usepackage{lineno}

\begin{document}
\onecolumn  

\title{Dynamical Systems with  a Cyclic Sign Variation Diminishing Property\thanks{Research  supported in part 
by  research grants from  the Israel Science Foundation  and the 
 US-Israel Binational Science Foundation.  
}}

 \author{Tsuff Ben Avraham, Guy Sharon,   Yoram Zarai and  Michael Margaliot\thanks{
The authors are  with the School of Electrical  Engineering,
Tel-Aviv University, Tel-Aviv~69978, Israel.
Corresponding author: Michael Margaliot (e-mail: \texttt{michaelm@eng.tau.ac.il})
 }}

\maketitle

\begin{abstract}
Several studies analyzed  certain nonlinear dynamical systems by showing that  
  the cyclic number of sign variations in the vector of derivatives is
	an integer-valued
	Lyapunov function.
	These results are based on direct analysis of the   
	dynamical equation satisfied by the vector of derivatives, i.e. the variational system. 
	However, it is natural to  assume  that
	they follow from the fact that the transition matrix in the   variational system satisfies
	a variation diminishing property~(VDP) with respect to the 
	cyclic number of sign variations in a vector. 
	Motivated by this, we develop the theoretical framework 
	of linear time-varying systems whose solution  satisfies 
	a~VDP with respect to the 
	cyclic number of sign variations. 
	This provides  an analogue  of the work of Schwarz 
		on \emph{totally positive differential systems}, i.e. linear time-varying systems whose solution  satisfies 
	a~VDP
	 with respect to the 
	standard (non-cyclic) number of sign variations.
\end{abstract}


\section{Introduction}
Let~$\M^+ \subset  \R^{n\times n}$ [$\M \subset  \R^{n\times n}$] denote the subset of~$n\times n$ real 
matrices that are tridiagonal with positive [nonnegative]  entries on the super- and sub-diagonals.  
In an interesting paper, Smillie~\cite{smillie} considered the  nonlinear system
\be\label{eq:fdot}
\dot y=f(y),
\ee
with~$y(t)\in \R^n$,    
satisfying that its Jacobian~$J(y):=\frac{\partial  f}{\partial y}(y) \in \M^+$ for all~$y$.
He showed that every trajectory of such a system either leaves any compact set or converges to an equilibrium. 
This result has found many applications as well as several interesting generalizations 
(see, e.g.~\cite{RFM_stability,chua_roska_1990,Donnell2009120,periodic_tridi_smith,fgwang2013}). 
 Smillie's analysis is based on showing that the number of sign variations in the vector of derivatives~$z(t):=\dot y(t)$ can only decrease with time. This was  
done by direct analysis of the differential equation for~$z(t)$, namely,~$\dot z=J(y)z$. 

Recently, it has been shown~\cite{fulltppaper}   that Smillie's results are intimately
 related to the pioneering, yet forgotten, work of Schwarz~\cite{schwarz1970} on
linear \emph{totally positive differential systems}~(TPDSs).\footnote{We use here a slightly different notion than in~\cite{schwarz1970} that
agrees with more modern terminology in the field 
 of totally positive   and totally nonnegative matrices.} 
Recall that a matrix is called totally positive~(TP)
[totally nonnegative~(TN)] if all its minors are positive [nonegative]. 
These matrices enjoy a rich and beautiful theory~\cite{total_book,pinkus}. 
In particular, multiplying a vector by a TP matrix can only decrease the number of sign variations in the vector. 

Schwarz considered  the linear time-varying 
system
\be\label{eq:ltv}
\dot x(t)=A(t)x(t),
\ee
with~$A(t)$ a continuous matrix function of~$t$. He called this system a~TPDS
 on a time interval~$(a,b)$ if its transition matrix~$\Phi(t,t_0)$ is totally positive
for any pair~$(t_0,t)$ with~$a<t_0<t<b$. Here the transition matrix is the matrix satisfying~$x(t)=\Phi(t,t_0)x(t_0)$.
In particular,~$\Phi(t_0,t_0)=I$. In the special case where~$A(t)$
is a constant  matrix, i.e.~$A(t)\equiv A$ then~$\Phi(t,t_0)=\exp((t-t_0)A)$. 
Of course, the transition matrix is real, square,  and 
nonsingular. 

Schwarz showed that if~$A(t)$ is a continuous matrix function of~$t$ then a necessary and sufficient condition for~TPDS is 
that~$A(t)\in \M$ for all~$t\in(a,b)$, and every  entry on the sub- or super-diagonal of~$A(t)$ is not zero on a time interval. In the particular case of a constant
 matrix~$A(t)\equiv A$
this means that~$\dot x(t)=Ax(t)$ is~TPDS  if and only if~(iff)~$A\in \M^+$.

Schwarz~\cite{schwarz1970} also showed that if~\eqref{eq:ltv} is~TPDS then the number of sign changes in~$x(t)$ can only decrease with time. 
To explain this, we recall three definitions for the number of sign variations  in a vector.
For a vector~$y\in\R^n$ with no zero entries, let~$\sigma(y)$ denote the number of indexes~$k\in\{1,\dots,n-1\}$ such that~$ y_k y_{k+1}<0$. 
By continuity, it is possible to extend the domain of definition of~$\sigma$ to the set
\begin{align}\label{eq:defvsret}
			V & :=\{y\in\R^n |    y_1\not =0, y_n\not =0, \\& \text{ if } y_i=0 \text{ for some } 2\leq i \leq n-1
			\text{ then } y_{i-1}y_{i+1}<0  \}. \nonumber 
\end{align}
For example~$y:=\begin{bmatrix} 1 &\varepsilon&  -1 \end{bmatrix}' \in V$ and~$\sigma(y)=1$ for all~$\varepsilon \in \R$ (including~$\varepsilon=0$). Two more definitions for the number of sign variations
in a vector, that are well-defined for any~$y\in \R^n$, are 
\begin{align*}
				s^-(y)&:=\sigma(\bar y), 
				\end{align*}
				 where $\bar y $   is the vector obtained from~$y$  by deleting all zero entries,
and
\begin{align*}
				s^+(y)&:=\max_{z\in P(y)}\{\sigma(  z)\},
\end{align*}
where~$P(y)$ is the set of vectors obtained 
 by replacing each zero entry in~$y$ by either $1$   or~$-1$. 
Clearly,~$s^-(y)\leq s^+(y)$ for all~$y$.
For example, for~$y=\begin{bmatrix}  0& 1& -2 \end{bmatrix}'$, $s^-(y)=1$ and~$s^+(y)=2$. Let
\[ W:=\{y\in\R^n|s^-(y)=s^+(y)\} .\]
 It is straightforward to show that~$W=V$. 

A classical and important   result from the theory of TP
 matrices~\cite{total_book}   states
	that if~$A$ is~TP  then
\[
	s^+(Ax)\leq s^-(x),\quad \text{for all } x\not =0.
\]
   At this point we can already see the connection between this
  sign variation diminishing property~(VDP) of~TP matrices
and the work of Smillie. Indeed, recall that~$z:=\dot y=f(y)$, and thus, 
\begin{align}\label{eq:zdot}
\dot z=\frac{\partial f}{\partial y}(y) \dot y=J(y)z.
\end{align}
This is     the \emph{variational system} associated with~\eqref{eq:fdot}. 
The assumptions of Smillie on the structure of~$J$ imply that~\eqref{eq:zdot} is a TPDS
and thus for any~$t>t_0$, 
\begin{align*}
s^+(z(t))&  = s^+(\Phi(t,t_0)z(t_0)) \\ &\leq s^-(z(t_0)).
\end{align*}

  From this it follows that~$\sigma(z(t)) \in V$ 
for all~$t$, except perhaps for up to~$n-1$ time points~$t_i$,
 and that at these points
\[
				\sigma(z(t_i^+))<\sigma(z(t_i^-)),
\]
see e.g.~\cite{schwarz1970,fulltppaper}.

As briefly mentioned
 in~\cite{fulltppaper}, several papers analyzed various properties of nonlinear dynamical systems
by showing that the number of \emph{cyclic} sign variations in the vector of derivatives~$z(t)$ is nonincreasing with time (see, 
e.g.,~\cite{Fusco1990,smith_sign_changes,poin_cyclic}). 
These results were proven directly, 
but it is natural to speculate that they are related to 
the fact that the  transition matrix of the variational system satisfies
 a cyclic variation diminishing property~(CVDP). 
In this paper, we develop the theoretical framework of such systems. 
  We say that the 
linear time-varying system~$\dot x(t)=A(t)x(t)$ is a  
\emph{cyclic variation diminishing  differential system}~({\model})  
if its transition matrix cannot increase the number of \emph{cyclic} sign variations. 
This is the ``cyclic analogue'' of a~TPDS.

The first step is to address the following question:
when does multiplication by a matrix~$A$ 
cannot increase  the  cyclic number of sign variations in a vector? 
Schoenberg and Whitney~\cite{schoenberg51} already addressed this question and subsequent work of 
Karlin~\cite[Ch.~5]{karlin_tp}  
    includes important characterizations of such matrices (and, more generally, kernels).
		However, the emphasis is on matrices~$A \in \R^{n\times m}$, with~$n>m$,
			whereas for the case of dynamical systems the relevant case is     square and nonsingular transition matrices. 
			We provide a simple   necessary and sufficient condition for a nonsingular square matrix~$A$ to satisfy
		a~CVDP (see Thm.~\ref{thm:lmainsccp} below).  
		
		The next step is to
		consider the matrix differential equation $\dot \Phi(t)=A(t)\Phi(t)$.
		For a constant  matrix~$A$, we provide a necessary and sufficient condition for 
		the transition matrix~$\exp((t-t_0)A)$, with~$t>t_0$,
			to satisfy  a~CVDP    (Thm.~\ref{thm:mecsauc}). 
			We then extend this to the case where~$A(t)$ is continuous in~$t$ 
			(Thm.~\ref{thm:nsct}). In the more general case 
		where~$t\to A(t)$ is measurable (but not necessarily continuous)  we provide
		a sufficient condition for the  transition matrix
			to satisfy  a~CVDP  (Thm.~\ref{thm:meas}).
			We also describe the implications of~{\model} to the solution of the vector equation~$\dot x=Ax$ (Thm.~\ref{thm:exp_is_ctp}).

		The remainder of this paper is organized as follows. The  next section reviews known 
		definitions and results  from the theory of~TP  matrices that will be used later on. 
		Section~\ref{sec:main} describes our main results. 
		We demonstrate one application of these results to
	  a nonlinear model called the \emph{ribosome flow model on a ring}. 
		Section~\ref{sec:conc} concludes and discusses directions for future research.

We use standard notation. Vectors [matrices] are denoted by small [capital] letters. 
$\R^n$ is the set of vectors with~$n$ real coordinates.
For a (column) vector $x\in\R^n$, $x_i$ is the $i$th entry of $x$, and~$x'$ is the transpose of~$x$.
For a matrix~$A$, $\tr(A)$ denotes  the trace of~$A$. 
A square matrix~$B$ is called Metzler if every off-diagonal entry of~$B$ is nonnegative. 
 The square identity matrix is denoted by~$I$, with dimension that should be clear from context.

	\section{Preliminaries}\label{sec:preli}
		We briefly review known definitions and results from  the rich and beautiful
		theory of totally nonnegative and totally positive 
		matrices that will be used later on.  For more information and proofs
		we refer to the   excellent monographs~\cite{total_book,pinkus,gk_book}.
		Unfortunately, this field suffers from nonuniform terminology. We follow the more modern terminology as in~\cite{total_book}.

We begin with some notation for the minors of a matrix~$A\in \R^{n\times m}$. 
Pick~$r \in\{1,\dots,\min\{n,m\}\}$, and let~$\alpha$ [$\beta$]
denote a set of~$r$  
integers~$1\leq i_1<\dots<i_{r}\leq n$
[$1\leq j_1<\dots<j_{r}\leq m$]. Then the minor of~$A$ corresponding to the rows indexed by~$\alpha$
and columns indexed by~$\beta$ is denoted~$A(\alpha|\beta)$. For example, for~$A=\begin{bmatrix} 1&2&3\\4&5&6\\7&8&9\end{bmatrix}$,
 $\alpha=\{1,3\}$
and~$\beta=\{1,2\}$, $A(\alpha|\beta)=\det( \begin{bmatrix} 1& 2  \\ 7&8 \end{bmatrix})= -6$. 

Pick~$A\in\R^{n\times m}$,
$B\in \R^{m\times p}$, and let~$C:=AB$.
The Cauchy-Binet formula~\cite[Ch.~1]{total_book} asserts that for any two
sets~$\alpha\subseteq\{1,\dots,n\}$, $\beta \subseteq \{ 1,\dots,p\}$ 
with the same cardinality~$k\in\{1,\dots,\min\{ n,m,p\}\}$, 
\be\label{eq:cbfor}
				 C(\alpha|\beta)= \sum_{  |\gamma| =k}  A(\alpha|\gamma)  B(\gamma |\beta) .
\ee
Here the sum is over all~$\gamma=\{i_1,\dots,i_k\}$, with~$1\leq i_1<\dots<i_k\leq m$.  
Thus, every minor of~$AB$ is the sum of products of minors of~$A$ and~$B$.
For example, for~$n=m=p$ and~$k=n$ Eq.~\eqref{eq:cbfor}
 gives the well-known formula~$\det(AB)=\det(A)\det(B)$.

We now turn to review the VDPs of~TP matrices or, more generally, of sign-regular matrices.
\begin{Definition}\label{def:ssr}
A matrix~$A\in\R^{n\times m}$
is called \emph{sign-regular of order~$k$} (denoted~$SR_k$)
if all minors of order~$k$  have the same \emph{non-strict} sign.
It  is called  \emph{strictly sign-regular of order~$k$}~($SSR_k$) 
if all  its minors of order~$k$  
are non-zero and have the same sign.
		It
		is called  \emph{strictly sign-regular}~(SSR) if it is~$SSR_k$ for all~$k\in \{1,\dots,\min\{n,m\}\}$, that is, 
		 all its minors of a given size are non-zero and share a common sign (that may vary from size to size).  
\end{Definition}
For example, the matrix~$\begin{bmatrix} 1&2\\ 0& 0\end{bmatrix}$ is~$SR_1$ and~$SR_2$,
and~$\begin{bmatrix} 1&2\\ 3& 1\end{bmatrix} $ is~SSR.

	Clearly, a TP matrix is~SSR. A classical result (see, e.g.~\cite[Ch.~V]{gk_book})    states
	that  a matrix~$A\in \R^{n\times m}$, with~$n>m$, 
	is~SSR  iff it satisfies the \emph{strong} sign variation diminishing property~(SVDP): 
	\[
					s^+(Ax)\leq s^-(x), \quad \text {for all } x\in\R^m\setminus \{0\}.
\]
Similarly, a square and \emph{nonsingular} matrix~$A\in \R^{n\times n}$
	is~SSR iff it satisfies the~SVDP. Note that this is not true for singular matrices.
	For example, it is straightforward to verify that~$A=\begin{bmatrix} 2 &2 \\1&1\end{bmatrix}$ satisfies
	the~SVDP, but it is not~SSR.
	

When using sign-regular  matrices to study dynamical systems, it is important to bear in mind that  in general 
		    the signs of minors are not invariant under similarity  transformations. 
		An important exception  is positive diagonal scaling. Indeed, if~$D $ is a diagonal matrix with positive
		entries on the diagonal then multiplying a matrix~$A$ by~$D$ either on the left or
right   changes the sign of no minor of~$A$, so in particular~$DAD^{-1}$ is SR [SSR] if and only if~$A$ is~SR [SSR].

Our first goal is to characterize square and nonsingular 
matrices that satisfy  an~SVDP with respect to
the \emph{cyclic} number of sign variations.

\subsection{Cyclic number of sign variations}
	For~$y\in\R^n$, let
	\begin{align}\label{eq:defsc-}
	s^-_c(y)& :=\max_{ i\in\{1,\dots,n\} } 
	s^-(\begin{bmatrix}  y_i & \dots &y_n &y_1& \dots&y_i
	\end{bmatrix}'),
	\end{align}
	This can be explained as follows: place the entries of~$y$ along a circular ring so that~$y_n$ is followed by~$y_1$, 
	then count~$s^-$ starting  from 
	any entry along the ring, and find the maximal value. 
	
	For example, for~$y=\begin{bmatrix} 0& 1 &0 & -3 \end{bmatrix}'$,
	$s^-_c(y)=s^-(\begin{bmatrix}    1 &0& -3& 0& 1  \end{bmatrix}') = 2$. 
	Similarly,
	\begin{align*}
	 s^+_c(y) & :=\max_{ i\in\{1,\dots,n\} }  
				s^+(\begin{bmatrix}  y_i & \dots &y_n &y_1& \dots&y_i
	\end{bmatrix}'),
	\end{align*}
	but here if~$y_i=0$ then in the calculation of~$s^+(\begin{bmatrix}  y_i & \dots &y_n &y_1& \dots&y_i
	\end{bmatrix}')$ \emph{both}~$y_i$s are replaced by either~$1$ or~$-1$. 
	For example, for~$y=\begin{bmatrix} 0& 1 &0 & -3 \end{bmatrix}'$,~$s^+_c(y)=s^+(\begin{bmatrix}    1&0&-3&0&1 \end{bmatrix}') =2$.
	Note that~$s^-_c(y) \leq s^+_c(y)$ for all~$y\in\R^n$,
	and that both~$s^-_c(y),s^+_c(y)$ are invariant 
	under  cyclic shifts of the vector~$y$.
	
	There is a simple and useful
	relation between the non-cyclic and cyclic number of sign variations of a vector. 
\begin{Lemma}\label{lem:rela}
For any vector~$x$,
\[
				s^-_c(x)=\begin{cases} 
							s^-(x), & \text{ if }s^-(x) \text{ is even},\\
							s^-(x)+1, & \text{ if }s^-(x) \text{ is odd}, 
							\end{cases}
\]
and, similarly, 
\[
				s^+_c(x)=\begin{cases} 
							s^+(x), & \text{ if } s^+(x) \text{ is even},\\
							s^+(x)+1, & \text{ if } s^+(x) \text{ is odd}.
									\end{cases}
\]
\end{Lemma}
For the sake of completeness, we include a proof of this result.

{\sl Proof of Lemma~\ref{lem:rela}.} 
Pick~$x\in\R^n$. If~$x=0$ then clearly~$s^-(x)=s^-_c(x)=0$, and~$s^+(x)=n-1$, furthermore,  
if~$n$ is even then~$s^+_c(x)=n$, and if~$n$ is odd then~$s^+_c(x)=n-1$. Thus, in this case
Lemma~\ref{lem:rela} holds.

Now  consider the case where~$x \in \R^n\setminus\{0\}$. 
Let~$p:=s^-(x)$. We may assume that the first non-zero entry of~$x$ is positive. 
Then the entries of~$x$ can be divided into~$p+1$ groups:
$
(x_1,\dots,x_{v_1})$, $(x_{v_1+1},x_{v_1+2},\dots,x_{v_2})$, $\dots$,
$(x_{v_p+1},x_{v_p+2},\dots,x_{v_{p+1}})$,
 where~$x_1,\dots,x_{v_1}\geq 0$ (with at least one of these entries positive),
$x_{v_1+1}<0$, $x_{v_1+2},\dots, x_{v_2}\leq 0$, $x_{v_2+1}>0$, and so on.
If~$p$ is even then~$x_{v_p+1 } >  0$  and~$x_{v_p+2},\dots,x_{v_{p+1}}\geq 0$.
Thus, the signs of the first and last group agree and~\eqref{eq:defsc-} yields~$s^-_c(x)=s^-(x)$.
If~$p$ is odd then~$x_{v_p+1} <  0$  and~$x_{v_p+2},\dots,x_{v_{p+1}}\leq 0$, so~$s^-_c(x)=s^-(x)+1$. The proof for~$s^+_c(x)$ is similar.~\hfill{$\square$}

Note that  Lemma~\ref{lem:rela}
 implies in particular that~$s_c^-(x),s^+_c(x)$ is always an even number. 
Furthermore, for~$x\in\R^n$, $s^-(x),s^+(x)$ take values in~$\{0,1,\dots,n-1\}$,
so we conclude that
\be\label{eq:nodeeevem}
s^-_c(x) ,s^+_c(x) \in \begin{cases} \{0,2,4,\dots,n \}, & \text{ if } n \text{ is even},\\
																		\{0,2,4,\dots,n-1 \}, & \text{ if } n \text{ is odd}.
\end{cases} 
\ee

We now consider the relation between non-cyclic and cyclic~VDPs. 
Suppose that~$A\in \R^{n \times m}$ and~$x\in \R^m $ satisfy
\be\label{eq:forsome}
s^+(Ax)\leq s^-(x) .
\ee
  If~$s^+(Ax)$ is even
 then~$s^+_c(Ax)= s^+(Ax)\leq s^-(x) \leq s^-_c(x) $. 
If~$s^+(Ax)$ is odd and~$s^-(x)$ is odd 
 then~$s^+_c(Ax)= s^+(Ax)+1 \leq s^-(x) +1 = s^-_c(x) $.
If~$s^+(Ax)$ is odd and~$s^-(x)$ is even  then~\eqref{eq:forsome} gives~$s^+(Ax)+1 \leq s^-(x)$, 
so~$s^+_c(Ax)= s^+(Ax)+1 \leq s^-(x)  = s^-_c(x) $.
In all  cases we see that~\eqref{eq:forsome} implies that
\[
s^+_c(Ax)\leq s^-_c(x). 
\]
We conclude that  if~$A$ satisfies a   VDP with respect to (w.r.t.) a non-cyclic number of sign variations it
also satisfies the same~VDP w.r.t. the cyclic number. 

However, it turns out that a weaker property is enough. 
To show this consider  a matrix~$A\in \R^{n\times n}$ that is~SSR. Then~$
			s^+(Ax)\leq s^-(x)$ for all~$ x\in \R^n \setminus\{0\}$,
		i.e.~$A$ satisfies the~SVDP.
As noted above, this  means that~$A$ also satisfies the strong cyclic VDP~(SCVDP):
\[
			s^+_c(Ax)\leq s^-_c(x),\quad\text{for all } x\in \R^n\setminus\{0\}.
\]
Let~$P_1,P_2\in\R^{n\times n}$ be cyclic permutation matrices. 
Since~$s^+_c(x),s^-_c(x)$ are invariant under cyclic permutations of~$x$,
we have that for all~$z\in \R^n\setminus\{0\}$,
\begin{align*}
s^+_c(P_1 A P_2^{-1} P_2 z)&=s^+_c( A  z)\\
& \leq  s^-_c(    z)\\
&= s^-_c( P_2  z).
\end{align*}
Thus, we conclude that~$B:=P_1 A P_2^{-1}$ also satisfies the~SCVDP, but~$B$ does not necessarily satisfy the~SVDP.
The next example demonstrates   this. 

\begin{Example}\label{exa:abexa}
For~$n=3$, consider the matrix 
$A:=	\begin{bmatrix}
		 5 &4 &1 \\
     4 &6 &4 \\
     1 &4 &5
	\end{bmatrix}.
$
This matrix is TP and thus satisfies   both the~SVDP and the~SCVDP.
Let~$P_1:=\begin{bmatrix}
		 0 &1 &0 \\
     0 &0 &1 \\
     1 &0 &0
	\end{bmatrix} 
	$ and~$P_2:=I$. 
	Then~$B:=P_1AP_2^{-1} = \begin{bmatrix}
		 4 &6 &4 \\
     1 &4 &5 \\
     5 &4 &1
	\end{bmatrix}$. This matrix satisfies the~SCVDP. However, 
	$B$ is not SSR (it has both positive and negative minors of order~$2$), and thus it does not satisfy the~SVDP.
	\end{Example}

Before ending this section, we state a well-known and  important  result that will be used later on. 

\begin{Proposition}\label{prop:sv1}
Consider a set of~$m$ vectors~$u^1,\dots,u^m \in \R^n$, with~$m < n$. 
Define the matrix~$U\in \R^{n\times m }$ by
\[
U:=\begin{bmatrix} u^1& u^2&\dots&u^m \end{bmatrix}.
\]
The following  two conditions are equivalent:
\begin{enumerate}[(1)]
\item	\label{cond:cis} For any~$c_1,\dots,c_m\in \R$, that are not all zero,
\be\label{eq:suim}
					s^+(\sum_{k=1}^m c_i u^i)\leq m-1.
\ee
\item			\label{cond:mino}
The matrix~$U$ is $SSR_m$.
\end{enumerate}
\end{Proposition}
For  a proof, see e.g.,~\cite{fulltppaper}.

\begin{Remark}\label{rem:overcon}
Note that the assumption that~$m<n$ cannot be dropped. 
For example, if~$m=n$ then condition~\eqref{cond:cis} 
always holds, whereas  condition~\eqref{cond:mino} holds iff~$U$ is nonsingular.
\end{Remark}

\section{Main Results}\label{sec:main}
Our first goal is to provide a necessary and sufficient condition for a square nonsingular matrix to satisfy the~SCVDP. 
We begin by stating an auxiliary result that will be used later on.
\subsection{Non-standard   VDP} 
We derive a   necessary and sufficient condition 
for a \emph{square and  nonsingular} matrix to satisfy       a non-standard~VDP. 
This result seems to be new
and    may be of independent interest, as it gives for any value~$r$
a clear interpretation of the~$SSR_r$ property   in terms of  this non-standard~VDP.  
\begin{Theorem}\label{thm:linspecp}
Let~$A\in\R^{n\times n}$ be a nonsingular matrix. 
  Pick~$p \in \{0,\dots,n-1\}$. 
Then the following  two conditions are equivalent:
\begin{enumerate}[(1)]
\item	\label{cond:onedip} For any vector~$c\in\R^n\setminus\{0\}$ with~$ s^-(c) \leq p$,
\be\label{eq:spsmalspt}
					s^+(Ac)\leq p.
\ee
\item			\label{cond:secdip}
$
A 
$
is~$SSR_{p+1}$. 
\end{enumerate}
\end{Theorem}

\begin{Example}\label{exa:proof1}
Consider the case~$p=0$. 
If~$c\in\R^n\setminus\{0\}$ satisfies~$ s^-(c) \leq 0$, i.e.~$s^-(c)=0$ then we may assume that every entry of~$c$ is nonnegative, and 
at   least one 
entry  is positive (as~$c\not =0$).  Then~$Ac$ is a nonnegative combination of the
columns of~$A$ with at least 
 one column taken with a positive weight. Also,~$Ac\not =0$, as~$A$ is nonsingular.

If condition~\eqref{cond:secdip} holds i.e.~$A$ is~$SSR_1$ then we may assume that
all  entries of~$A$ are positive. Thus, every entry of~$Ac$ is positive, so~$s^+(Ac)=0$
and condition~\eqref{cond:onedip} holds.

Suppose that  condition~\eqref{cond:onedip} holds. Taking~$c=e^k$, where~$e^k$ is the $k$th canonical vector in~$\R^n$, we
get that~$s^+(A e^k)=0$, i.e. all the entries in column~$k$ of~$A$ have the same strict sign. Let~$\varepsilon_k \in\{-1,1\}$ denote this sign.
 Seeking a contradiction, assume that there exist~$i,j$ such that~$\varepsilon_i=1$ and~$\varepsilon_j=-1$. 
Then we can find~$d_i,d_j>  0$ such that for~$c=d_i e^i+d_j e^j$ the vector~$Ac=d_iAe^i+d_jAe^j$ includes a zero entry, so~$s^+(Ac) >0$. 
This contradicts condition~\eqref{cond:onedip}, so we conclude that all the~$\varepsilon_k $s are identical, i.e.~$A$ is~$SSR_1$.

\end{Example}

\begin{Example}\label{exa:pott}
Consider the 
matrix
\be\label{eq:ant}
A=\begin{bmatrix}  1& 2 & 0 &0 \\
											0&1&1&0 \\ 0 & 0&2 &1\\1&0&0 &2
\end{bmatrix}.
\ee
 It is straightforward to verify that this matrix is nonsingular and~$SSR_3$,
i.e. Condition~\eqref{cond:secdip}
  holds for~$p=2$.  
We will show that Condition~\eqref{cond:onedip}   holds for this value of~$p$.
Pick  a vector~$c\in\R^4\setminus\{0\}$ with~$ s^-(c) \leq 2$. 
  Note that
\begin{align*}
Ac= \begin{bmatrix}  c_1+2c_2&c_2+c_3& 2 c_3+c_4&c_1+2c_4   \end{bmatrix}'.
\end{align*}
Seeking a contradiction, assume that~$s^+(Ac)>2$, i.e.~$s^+(Ac)=3$. 
Then without loss of generality we can assume that
\begin{align}\label{eq:fghtyyr}
  c_1+2c_2 &\geq 0, \nonumber \\
	c_2+ c_3 &\leq 0,\\
	2 c_3+c_4&\geq 0,\nonumber\\
	c_1+2 c_4 &\leq 0\nonumber.   
\end{align}
We consider three cases.

\noindent{\sl Case 1.} Suppose that~$c_1< 0$. Then the first equation in~\eqref{eq:fghtyyr}
yields~$c_2> 0$, the second equation gives~$c_3<0$, and  the third yields~$c_4> 0$. 
But this means that~$s^-(c)=3$ and this is a contradiction. 

\noindent{\sl Case 2.} Suppose that~$c_1= 0$. Then the first three equations in~\eqref{eq:fghtyyr}
yield~$c_2\geq  0$, $c_3\leq 0$, and~$c_4\geq  0$. Now the fourth equation gives~$c_4=0$. 
Substituting this in the third equation gives~$c_3=0$, and the second equation gives~$c_2=0$.
 We conclude that~$c=0$, which is a contradiction. 

 \noindent{\sl Case 3.} Suppose that~$c_1> 0$. Then the fourth   equation in~\eqref{eq:fghtyyr}
yields~$c_4 <   0$, the third gives~$c_3> 0$, and the second gives~$c_2<0$.
This means that~$s^-(c)=3$, which is a contradiction. 

Summarizing, in each case we obtained a contradiction to~$s^+(Ac)=3$, so~$s^+(Ac)\leq 2$. 
\end{Example}

 We emphasize that condition~\eqref{cond:onedip} in Thm.~\ref{thm:linspecp} does not  assert  that~$s^-(c)\leq p$
implies that~$s^+(Ac)\leq s^-(c)$,
but only that~$s^+(Ac)\leq p$.  
For example, for the matrix~$A$ in~\eqref{eq:ant} and the vector~$c=\begin{bmatrix}
0&0&0&-1\end{bmatrix}'$ we have~$s^-(c)=0\leq  2 $, and
 $s^+(Ac)=  s^+(\begin{bmatrix} 0&0&-1&-2\end{bmatrix}') = 2 $.

{\sl Proof of Thm.~\ref{thm:linspecp}.}
First note that for~$p=n-1$ both conditions~\eqref{cond:onedip} and~\eqref{cond:secdip}
in Thm.~\ref{thm:linspecp} 
  hold, as~$s^+(z)\leq n-1$ for all~$z\in\R^n$,
and~$A$ is nonsingular and thus~$SSR_n$. Also, we already proved Thm.~\ref{thm:linspecp}
 when~$p=0$. 
Thus, we need to prove the result for any~$p\in\{1,\dots,n-2\}$.

Assume that condition~\eqref{cond:onedip} holds for  some value~$p$ in the
 range~$ \{1,\dots,n-2\}$.
We will show that all minors  of  order~$p+1$ of~$A$  are non-zero and have the same sign.
 Pick  a set of~$p+2$ indices~$1\leq k_1<k_2<\dots <k_{p+2}\leq n$. 
For any~$c\in\R^n \setminus \{0\}$ with~$s^-(c)\leq p$
define~$\bar c\in \R^n$
by:
\[
				\bar c_i:=\begin{cases} c_i, & \text{if }  i\in\{k_1,\dots,k_{p+1}\}, \\
				                          0 ,& \text{otherwise}. 
																	\end{cases} 
\] 
Then~$s^-(\bar c)\leq p$, as~$\bar c$ has no more than~$p+1$ non-zero entries. 
Let~$a^i\in\R^n$ denote the~$i$th column of~$A$.  Then~$A \bar c =\sum_{i=1}^n \bar{c}_i a^i =\sum_{j=1}^{p+1} c_{k_j}a^{k_j}$. 
Applying condition~\eqref{cond:onedip} 
 implies that for any~$\bar c \not = 0$, $s^+(A \bar c)\leq p  $. 
This means that the set~$\{a^{k_1},\dots,a^{k_{p+1}}\}$ 
satisfies condition~\eqref{cond:cis}  in  
Prop.~\ref{prop:sv1} (note that~$n>p+1$). Thus, all minors of the form
\be\label{eq:aaminors}
A( i_1\;\dots\; i_{p+1}| k_1\; \dots \;k_{p+1}),  
\ee
 with $1\leq i_1<i_2<\dots<i_{p+1}\leq n$,
 are non-zero and have the same sign. Denote this sign by~$\varepsilon(k_1, \dots ,k_{p+1})\in\{-1,1\}$.
 It remains to show that this sign
depends  on the value~$p$, but not on the particular choice of~$k_1,\dots,k_{p+1}$.
Pick~$v\in\{1,\dots,p+1\}$.
We will show that
\begin{align}\label{eq:vere}
										&\varepsilon(k_1,\dots,k_{v-1},k_{v+1},\dots,k_{p+2}) \nonumber \\&= 
										\varepsilon(k_1,\dots,k_{v},k_{v+2},\dots,k_{p+2 }).
\end{align}
To do this, define~$p+1$ vectors~$ \{ \bar a^{k_1},\dots,\bar a^{k_v},\bar a^{k_{v+2}}
,\dots,\bar a^{k_{p+2}} \}  $ 
by
\[
					\bar a^{k_i}:=\begin{cases}  a^{k_i}, & i\not =v,\\
																			 d_v a^{k_v}+d_{v+1} a^{k_{v+1}}, &i=v,
																			\end{cases}
\]
 where~$d_v,d_{v+1}>0$.
Pick scalars~$\bar c_1,\dots,\bar c_v,\bar c_{v+2},\dots,\bar c_{p+2}$,
 that are not all zero, and 
let
\be\label{eq:deftra}
a:=\sum_ {i=1 \atop i\not = v+1}^{p+2} \bar  c_i \bar a^{k_i}. 
\ee
Then 
\be\label{eq:aprt}
				a=\sum_{i=1 }^{p+2} g_k  a^{k_i},
\ee
with~$g_v:=\bar c_v d_v$, $g_{v+1}:=\bar  c_{v} d_{v+1}$, and~$g_i=\bar c_i$ for
 all other~$i$. Let~$g:=\begin{bmatrix} g_1&\dots&g_{p+2}\end{bmatrix}'$.
Note that~$g\not =0$, and that since~$g_v g_{v+1}=(\bar c_v)^2 d_vd_{v+1}\geq 0$, 
$s^-(g)\leq p$. 
Applying condition~\eqref{cond:onedip}  to~\eqref{eq:aprt}
yields~$s^+(a)\leq p$, that is,~$s^+(\sum_ {i=1 \atop i\not = v+1}^{p+2} \bar  c_i \bar a^{k_i})\leq p$.
 Let~$\bar A \in\R^{n\times (p+1)}$ be the matrix
\begin{align*} 
\bar A:&=\begin{bmatrix} \bar a^{k_1}&\dots&\bar a^{k_{v-1}}& \bar a^{k_v}&
\bar a^{k_{v+2}}&\dots&\bar a^{k_{p+2}} \end{bmatrix}.
\end{align*}
Applying Prop.~\ref{prop:sv1} to the set of~$p+1$ 
vectors~$\bar a^{k_1},\dots,\bar a^{k_v},\bar a^{k_{v+2}},\dots,\bar a^{k_{p+2}}\in\R^n$ we conclude that
 all minors
\begin{align*} 
\bar A(  i_1,\dots,i_{p+1}  |  k_1,\dots,k_v,k_{v+2},\dots,k_{p+2} )
\end{align*}
are non-zero. Thus, all minors of the form
\begin{align}\label{eq:twomin}
&d_v A(  i_1,\dots,i_{p+1}  |  k_1,\dots,k_v,k_{v+2},\dots,k_{p+2} ) \nonumber\\
&+d_{v+1} A(  i_1,\dots,i_{p+1}  |  k_1,\dots,k_{v-1},k_{v+1},\dots,k_{p+2} )
\end{align}
are non-zero. This holds for  all~$d_v,d_{v+1}>0$,
so  the two minors in~\eqref{eq:twomin} have the same sign. 
Since this is true for all~$v\in\{1,\dots,p+1\}$,
we conclude that the sign~$
	\varepsilon(k_1,\dots,k_{p+1}) 
$
does not change if we change any one of the indices~$k_i$, and thus it is independent of
the choice of~$k_1,\dots,k_{p+1}$. 
This completes the proof that condition~\eqref{cond:onedip} implies condition~\eqref{cond:secdip}.

To prove the converse implication, assume that  
$
A =\begin{bmatrix} a^1&  \dots & a^n \end{bmatrix} 
$
is $SSR_{p+1}$ for some~$p \in\{1,\dots,n-2\}$.   
 Pick~$c\in \R^n\setminus\{0\}$ such that~$ s^-(c)\leq p$. Let~$k:=s^-(c)$.
We may assume that the first non-zero entry of~$c$ is positive. 
Then we can decompose~$c$ into~$k+1$ groups:
\begin{align} \label{eq:gdecp} 
 &(c_{1},c_{2},\dots,c_{v_{1}}), (c_{v_1+1},c_{v_1+2},\dots,c_{v_{2}}),\dots,\nonumber\\& (c_{v_k+1},c_{v_k+2},\dots,c_{v_{k+1}}),
\end{align}
where
 $c_1,\dots,c_{v_1}\geq 0$ (with at least one of these entries positive);
$c_{v_1+1}<0$, $c_{v_1+2},\dots, c_{v_2}\leq 0$, $c_{v_2+1}>0$, and so on. 
Define vectors~$u^1,\dots, u^{k+1} \in \R^{n}$ by
\[
				u^1:=\sum_{j=1}^{v_1}|c_j|a^j, \; u^2:= \sum_{j=v_1+1}^{v_2}|c_j|a^j,\dots.
\]
Then
\begin{align}\label{eq:asqre}
Ac&=\sum_{j=1}^{n}c_j a^j \nonumber\\
&= u^1-u^2+u^3-\dots+(-1)^k u^{k+1}.
\end{align}
Note that every~$u^i$ is a non-negative and non-trivial sum of a consecutive set of~$a^k$s. 
We now consider two cases.

\noindent {Case 1.}  
 Suppose that~$k=p$. Let
$ U:=\begin{bmatrix} u^1 &  \dots & u^{k+1}\end{bmatrix} . 
$ Note that the dimensions of~$U$ are~$n \times(p+1)$ with~$n>p+1$. 
Since~$A$ is~$SSR_{p+1}$,
all the minors of order~$(p+1)$ of~$U$ are non-zero and have the same sign.  
 Applying Prop.~\ref{prop:sv1} to~\eqref{eq:asqre}
yields~$s^+(Ac) \leq p $.

\noindent {Case 2.}
 Suppose that~$k < p$. 
 We then split the~$k+1$ groups in~\eqref{eq:gdecp} and
generate~$p+1$ vectors~$w^1,\dots,w^{p+1}$ by modifying~$u^1,\dots,u^{k+1}$ as follows.
Suppose for concreteness 
  that the first group in~\eqref{eq:gdecp} contains more than one element and that~$c_1>0$.  If~$c_2>0$  we replace~$u^1$ by two vectors
$w^{1}:=|c_1| a^1$ and~$w^{2}:=\sum_{j=2}^{v_1}|c_j|a^j$, and let~$w^{3}:=u^2$, $w^4:=u^3,\dots,w^{k+2}:=u^{k+1}$.  
Note that~$Ac=  (w^1+w^2)-w^3+w^4-w^5-\dots+(-1)^k w^{k+2}$.
If~$c_2=0$ 
we replace~$u^1$ by two vectors
$w^{1}:=|c_1| a^1$ and~$w^{2}:=a^2$, and let~$w^{3}:=u^2$, $w^4:=u^3,\dots,w^{k+2}:=u^{k+1}$. 
Note that~$Ac= (w_1+0 *w_2)-w^3+w^4-w^5-\dots+(-1)^k w^{k+2}$.
We can continue this decomposition, if necessary, using the entry~$c_3$, etc. 
If the first group in~\eqref{eq:gdecp} contains   only one element then we can follow the same idea using 
any entry in any group that contains more than a single element. In this way
  we can   generate a set~$w^1,\dots,w^{p+1}$
such that
\be \label{eq:sdcsp}
Ac=\sum_{j=1}^{p+1} d_i w^i,
\ee
 where the~$d_i$s are not all zero, and the matrix~$W:= \begin{bmatrix} w^1 &  \dots & w^{p+1}\end{bmatrix} \in \R^{n \times  (p+1) } $  
is~$SSR_{p+1}$ because~$A$ is~$SSR_{p+1}$.  
 Applying Prop.~\ref{prop:sv1} to~\eqref{eq:sdcsp}
yields~$s^+(Ac) \leq p $. 

Summarizing we showed  in both cases  that condition~\eqref{cond:secdip}
implies that condition~\eqref{cond:onedip} holds.
This completes the proof of Thm.~\ref{thm:linspecp}.~\hfill{$\square$}

\subsection{Conditions for cyclic VDP} 
The main  result in this subsection
  provides a necessary and sufficient condition for a
square  and nonsingular matrix to satisfy the~SCVDP. 

\begin{Theorem} \label{thm:lmainsccp}
Let~$A\in\R^{n\times n}$ be a nonsingular matrix.
The following  two conditions are equivalent:
\begin{enumerate}[(1)]
\item	\label{cond:scvggp} For any vector~$x\in\R^n\setminus\{0\}$,  
\be\label{eq:mainc1}
					s^+_c(Ax)\leq s^-_c(x).
\ee
\item			\label{cond:mainytr}
The matrix~$A$ is $SSR_r$ for all odd~$r$ in the range~$r\in\{1, \dots,n\}$.
\end{enumerate}
\end{Theorem}

\begin{Example}\label{exa:sopim}
Suppose that~$A\in\R^{3\times 3}$ is nonsingular. Then~$A$ is~$SSR_3$.
 Thm.~\ref{thm:lmainsccp} asserts that~$A$ satisfies the~SCVDP   if and only if
all the entries of~$A$ are either all positive or all negative. 
Note that this agrees  
with the results in Example~\ref{exa:abexa}.
\end{Example}

 {\sl Proof of Theorem~\ref{thm:lmainsccp}.}
Suppose that condition~\eqref{cond:mainytr} holds. Pick~$x\in\R^n\setminus\{0\}$.
Let~$k$ be such that~$s^-_c(x)=2k$ (recall that~$s^-_c$ and~$s^+_c$ are always even). 
If~$2k=n$ then clearly~\eqref{eq:mainc1} holds, so we may assume that~$2k\leq n-1$.. 
Since~$s^-_c(x)=2k$,~$s^-(x)\leq 2k$.
Condition~\eqref{cond:mainytr} implies in particular that~$A$ is~$SSR_{2k+1}$, and 
Thm.~\ref{thm:linspecp} yields~$s^+(Ax)\leq 2k$. Hence~$s^+_c(Ax)\leq 2k=s^-_c(x)$. 
This shows that condition~\eqref{cond:mainytr} implies condition~\eqref{cond:scvggp}.

To prove the converse implication assume that  condition~\eqref{cond:scvggp} holds. 
Pick an \emph{even} number~$p \in \{0,\dots,n-1\}$,
 and a vector~$x\in\R^n \setminus\{0\}$ such that~$s^-(x)\leq p$. 
Condition~\eqref{cond:scvggp}  yields~$s^+_c(Ax)\leq p$,   
so~$s^+(Ax)\leq p$. Thm.~\ref{thm:linspecp} implies that~$A$ is $SSR_{p+1}$.  
Since~$p$ is an arbitrary even number in~$\{0,\dots,n-1\}$, we conclude that
condition~\eqref{cond:mainytr} holds.~\hfill{$\square$}

Using a standard  continuity argument yields the following
condition for \emph{weak} cyclic~VDP.

\begin{Theorem}\label{thm:lwalajf}
Let~$A\in\R^{n\times n}$ be a nonsingular matrix.
The following  two conditions are equivalent:
\begin{enumerate}[(1)]
\item	\label{cond:wescvggp} For any vector~$x\in\R^n\setminus\{0\}$,
\be\label{eq:wamainc1}
					s^-_c(Ax)\leq s^-_c(x).
\ee
\item			\label{cond:wemainytr}
The matrix~$A$ is $SR_r$
for all odd~$r$ in the range~$r\in\{1, \dots,n\}$.
%
%
\end{enumerate}
\end{Theorem}

{\sl Proof of Thm.~\ref{thm:lwalajf}.}
Suppose that condition~\eqref{cond:wemainytr} holds.
For~$y \in\R$, let~$F(y)$  denote the~$n\times n$ matrix
whose~$i,j$ entry is~$\exp(- (i-j)^2 y) $. For example, for~$n=3$, 
\[
F(y)=\begin{bmatrix} 
                             1& \exp(-   y)& \exp(- 4  y)\\
														  \exp(-   y)&1& \exp(-    y)\\
														 \exp(-4   y)& \exp(-    y)&1
\end{bmatrix} .
\]
It is well-known that~$F(y)$ is~TP for all~$y>0$~\cite[Ch.~II]{gk_book},
and~$\lim_{y\to\infty}F(y)=I$.
Fix~$y>0$ and let~$F:=F(y)$,~$B:=F A$.
Pick an odd index~$r\in\{1,\dots,n\}$, and let~$\alpha,\beta$
denote two sets of~$r$ integers~$1\leq i_1<\dots<i_r\leq n$ and~$1\leq j_1<\dots<j_r\leq n$, respectively.
The Cauchy-Binet formula yields
\[
B(\alpha|\beta) =\sum_{\gamma} F (\alpha|\gamma)  A(\gamma | \beta),
\]
where the sum is over all~$\gamma=\{k_1,\dots,k_r\}$, with~$ 1\leq k_1<\dots< k_r\leq n $.
Using the fact that~$F$ is TP, all minors of order~$r$ of~$A$
have the same non-strict sign and they are not all zero (as~$A$ is nonsingular), we conclude
that~$B$ is~$SSR_r$.
Now Thm.~\ref{thm:lmainsccp} implies that for any~$x\in\R^n\setminus\{0\}$, $s^+_c( F(y) A x)\leq s^-_c (x)$,
so \[
s^-_c( F(y) A x)\leq s^-_c (x).
\]
Taking~$y \to \infty$  yields~\eqref{eq:wamainc1}.
   Thus, condition~\eqref{cond:wemainytr} implies condition~\eqref{cond:wescvggp}.

To prove the converse implication, assume that condition~\eqref{cond:wescvggp} holds.
Pick~$y>0$. Then~$F(y)$ is~TP and thus satisfies the~SVDP and thus the~SCVDP, so for any~$x\in\R^n \setminus\{0\}$,
\[
s^+_c(F(y) A x) \leq s^-_c(Ax).
\]
Combining this with~\eqref{eq:wamainc1} yields
\[
s^+_c(F(y) A x) \leq s^-_c(x).
\]
Now Thm.~\ref{thm:lmainsccp} implies
that~$F(y)A$ is~$SSR_r$ for all odd~$r$.
Taking~$y\to \infty$  we conclude that~$A$ is~$SR_r$ for all odd~$r$.~\hfill{$\square$}

\begin{Example}\label{exp:Ascminus}
Consider the case~$A\in\R^{2\times 2}$. In this case,
condition~\eqref{cond:wemainytr} in Thm.~\ref{thm:lwalajf} holds iff~$A$ is~$SR_1$ i.e. iff 
 all the entries in~$A$
have the same non-strict sign.  
We will show directly that in this case condition~\eqref{cond:wescvggp} in Thm.~\ref{thm:lwalajf} also holds. 
Pick~$x\in\R^2\setminus\{0\}$. We consider two cases (recall that~$s_c^-(x)$ is always even).

\noindent {\sl Case 1.} 
If~$s_c^-(x)=2$ then clearly~\eqref{eq:wamainc1} holds. 

\noindent {\sl Case 2.} 
Suppose that~$s_c^-(x)=0$. This means that~$x_1,x_2$ have the same non-strict sign.
 We may assume that~$x_1,x_2\geq 0$.
Then both entries of~$
Ax = \begin{bmatrix} a_{11}x_1+a_{12} x_2 \\
 a_{21}x_1+a_{22} x_2 \end{bmatrix}  
$ have the same non-strict sign,
so~$s_c^-(Ax)=0$ and again~\eqref{eq:wamainc1} holds. 
\end{Example}

We   now proceed to defining the analogue of a~TPDS for the case of the
 cyclic number of sign variations.

\subsection{Dynamical systems satisfying   a cyclic VDP} 
Consider the   linear time-varying system~\eqref{eq:ltv}. 
The associated matrix differential equation is
\be\label{eq:mdfe}
\dot  \Phi(t) =A(t)\Phi(t ).
\ee
Let~$\Phi(t,t_0)$ denote the solution of this equation
 satisfying~$\Phi(t_0, t_0)=I$.
\begin{Definition}
We say that~\eqref{eq:ltv} (and also~\eqref{eq:mdfe}) is a 
\emph{cyclic variation diminishing  differential system}~({\model}) on a time interval~$(a,b)$
if~$\Phi(t,t_0)$   satisfies the~SCVDP for any pair~$(t_0,t)$
with~$a<t_0<t<b$. 
\end{Definition}

\begin{Example}\label{exa:tbt}
Consider the fixed   matrix
$
			A = \begin{bmatrix} a_{11} & a_{12} \\a_{21} &a_{22} \end{bmatrix}$. 
			Then
\begin{align*}
	\Phi(t,t_0)&=		\exp( (t-t_0 )A )\\
	&=\begin{bmatrix}1+ a_{11}d+o(d) & a_{12}d+o(d)  \\a_{21}d+o(d) &1+a_{22}d+o(d) \end{bmatrix} ,
\end{align*}
where~$d:=t-t_0$. 
This implies that the following three statements  are equivalent: 
(1)~$a_{12},a_{21} >  0$; (2)~there exists~$\varepsilon>0$ such 
  that~$\exp((t-t_0)A)$ is~$SSR_1$ for all~$(t-t_0)\in (0,\varepsilon]$; and
	(3)~$\exp((t-t_0) A)$ is~$SSR_1$ for all~$t>t_0$.
	
	Combining this with Thm.~\ref{thm:lmainsccp}, we  conclude that 
	the system~$\dot x=Ax$, with~$A\in\R^{2\times2}$, is {\model} on any   
	time interval 
	  iff~$a_{12},a_{21}>0$. 
\end{Example}

The next result describes one implication  of~{\model}. 
Let
\be\label{eq:vcdef}
				V_c:=\{x\in\R^n: s^-_c(x)=s^+_c(x)\}.
\ee
	\begin{Theorem}\label{thm:exp_is_ctp}
Suppose that~\eqref{eq:ltv} is {\model} on~$(a,b)$.   
If~$x(t)$ is not the trivial solution~$x(t)\equiv  0$  then:
\begin{enumerate}[(1)]
\item  $s^-_c(x(t))$, $s^+_c(x(t))$ are   non-increasing functions of
 time on~$(a,b)$;
\item $x(t) \in V_c$  for all~$t\in (a,b)$, except perhaps for up  to~$\lfloor n/2 \rfloor $ discrete  values of~$t$; 
\item there exists a time~$T\geq 0$ such that
\[
				s^-_c(x(\tau))=s^+_c(x(\tau)),\quad \text{ for all } \tau \geq T,
\]
and   for all such~$\tau$ no two consecutive entries (in the cyclic order) of~$x$
satisfy
\[
x_i(\tau)=x_{i+1}(\tau)=0.
\]
\end{enumerate}
\end{Theorem}


{\sl Proof of Thm.~\ref{thm:exp_is_ctp}.}
  For any~$a<t_0<t<b$, we have~$x(t) =\Phi(t,t_0)x(t_0 )$.
	Since~$x(t_0)\not =0$ and~$\Phi(t,t_0)$    
	  satisfies the~SCVDP,
		\be\label{eq:slop}
	s^+_c(x(t))\leq s^-_c (x(t_0) ).
	\ee
	Thus,~$s^-_c(x(t))\leq s^+_c(x(t))\leq s^-_c (x(t_0) ) $.
	If~$x( t_0) \in V_c$ then~$s^-_c (x(t_0) )  = s^+_c (x(t_0) )$, so~\eqref{eq:slop} yields  
	\[
					s^+_c(x(t))  \leq  s^+_c (x(t_0) ) .
	\]
	If~$x(t_0)\not \in V_c$ then~$s^-_c (x(t_0) )  < s^+_c (x(t_0 ) )$, so 
	\be\label{eq:decft}
					s^+_c(x(t))  < s^+_c (x(t_0) ). 
	\ee
	Thus, $s^+_c (x(t))$ never increases, 
	and it strictly decreases as~$x(t)$ goes through a point that is not in~$V_c$.  
	Since~$s^+_c (x(t))$ is even the decrease  
	is by at least two.
	Since~$s^+_c$ takes  values  in~$\{0,1,\dots,n\}$ [$\{0,1,\dots,n-1\}$] when~$n$ is even [odd],
	 this implies
	that~$x(t)\in V_c$ for all~$t $, except perhaps for up to~$\lfloor n/2 \rfloor $  discrete  points.
	This proves the first two assertions in the theorem. The third assertion follows immediately 
	from the second and the fact that if~$y$ has
	two consecutive (in the cyclic order) zeros then~$s^-_c(y)< s^+_c(y)$.~\hfill{$\square$}

The  proof  of Thm.~\ref{thm:exp_is_ctp}
 shows that we may view~$s^- (x(t))$, $s^+(x(t))$ as     integer-valued Lyapunov functions of
 the time-varying linear system~\eqref{eq:ltv}. 


Using Thm.~\ref{thm:lmainsccp} yields  a   converse for
Thm.~\ref{thm:exp_is_ctp}.

\begin{Corollary}\label{coro:conv}
Suppose that the  solution of~\eqref{eq:ltv}   satisfies~\eqref{eq:slop} for 
 all~$a<t_0<t<b$ and all~$x(t_0)\in\R^n \setminus\{0\}$.
Then  for all~$(t_0,t)$ with~$a<  t_0 < t<b$ all the odd minors  of the transition matrix~$\Phi(t,t_0)$ are   positive. 
\end{Corollary}

{\sl Proof of Corollary~\ref{coro:conv}.}
Pick an odd~$r\in\{1,\dots,n\}$. The matrix~$\Phi(t_0,t_0)=I$ has a minor of order~$r$ that is one.
By Thm.~\ref{thm:lmainsccp}, all  minors of order~$r$ of~$\Phi(t,t_0)$, with~$t>t_0$,
are non-zero and have the same sign.   By continuity, this implies
that all minors of order~$r$ of~$\Phi(t,t_0)$, with~$t>t_0$,  are positive.
Since~$r$ is an arbitrary odd number, this completes the proof.~\hfill{$\square$}



The next natural question is what conditions on~$A(t)$ 
guarantee that~\eqref{eq:ltv} is~{\model}. The next subsection answers  this question 
for the case where~$A(t)\equiv A$. 
\subsection{The case~$A$ constant}
\begin{Definition}
Let~$\Q^+\subset \R^{n\times n}$ denote the set of~$n\times n$
matrices that   are Metzler, irreducible, and can have    non-zero entries   only
  on the main diagonal, super- and sub-diagonals and entries $(1,n)$
and~$(n,1)$.
\end{Definition}
 For example, for~$n=4$ the matrices:
\[
\begin{bmatrix}
									*& +& 0&0 \\ 0& *& +&0 \\ 0& 0& *&+ \\ + & 0& 0&*
\end{bmatrix},
\begin{bmatrix}
									*& +& 0&0 \\ +& *& +&0 \\ 0& +& *&+ \\ 0 & 0& +&*
\end{bmatrix},
\]
where~$*$ denotes some value and~$+$ denotes a positive value, 
are in~$\Q^+$. 

 Note that if~$A\in \Q^+$ then~$A'\in \Q^ +$, and that if~$A,B\in\Q^+$ then~$A+B \in \Q^+$. 

\begin{Theorem} \label{thm:mecsauc} 
Fix  a time interval~$(a,b)$. Let~$A$ be a constant~$n\times n$ matrix. 
The system~$\dot x=Ax$ is {\model}   on~$(a,b)$ 
iff~$A\in \Q^+$. 
\end{Theorem}

\begin{Example} 
Consider the case $n=2$. 
In Example~\ref{exa:tbt} we saw that 
	the system~$\dot x=Ax$, with~$A\in\R^{2\times2}$, is~{\model} on any   
	time interval 
	  iff~$a_{12},a_{21}>0$. On the other-hand, the definition of~$\Q^+$
		shows that~$A \in\Q^+$ iff~$a_{12},a_{21}>0$.
\end{Example}

Proving Theorem~\ref{thm:mecsauc}   requires introducing more notation. Before that we make several comments.
First, it follows from the definition of~$\Q^+$ that if~$A\in \M^+$ then~$A\in \Q^+$. 
This also makes sense, as~$A\in\M^+$ implies that~$\dot x=Ax$ is~TPDS and thus it is also~{\model}. 
If~$A\in (\Q^+ \setminus \M^+) $ then the system is not~TPDS and thus~$\exp(At)$ does not satisfy the~SVDP.
However, the next remark shows that any possible
 violation of the~SVDP has a particular structure. 
 \begin{Remark}\label{rem:QnotM}
Pick~$A\in (\Q^+ \setminus \M^+)$. 
Then~$\dot x = Ax$ is {\model},
 but not~TPDS. 
Suppose that at some time~$\tau$ there is an increase in~$s^-(x(t))$ (the analysis for~$s^+(x)$ is similar), say,
\be\label{eq:inctry}
s^-(x(\tau^-))=k \text{ and } s^-(x(\tau^+))> k.
\ee
Assume that~$k$ is even. Then  Lemma~\ref{lem:rela} yields
\begin{align*}
s^-_c(x(\tau^-))=k \text { and } s^-_c(x(\tau^+))> k ,
\end{align*}
but this is impossible as the system is~{\model}.
 We conclude that if~\eqref{eq:inctry} holds at some time~$\tau$ then~$k$
is odd, and by  Lemma~\ref{lem:rela},~$s^-_c(x(\tau^-))=k+1$.
Since the system is~{\model}, we must have~$s^-_c(x(\tau^+))\leq k+1$.
 Lemma~\ref{lem:rela}  now implies  that~$s^-(x(\tau^+))=k+1$. 
Summarizing,  if there is an increase in~$s^-(x(t))$ at time~$\tau$ 
then we must have
\begin{align}\label{eq:inctry2}
s^-(x(\tau^-))&=2i+1,\; & s^-(x(\tau^+))=2i+2,\nonumber \\
s^-_c(x(\tau^-))&=2i+2, & s^-_c(x(\tau^+))=2i+2,
\end{align}
for some~$i$. 
In particular, if for some time~$T$ we have~$s_c^{-}(x(t))=s_c^{+}(x(t))=0$ for all~$t\geq T$
 then~$s^{-}(x(t))=s^{+}(x(t))=0$ for all~$t\geq T$.
\end{Remark}
 
\begin{figure*}
\begin{subfigure}[t]{0.5\textwidth}
\begin{center}
\includegraphics[width= 8cm,height=7cm]{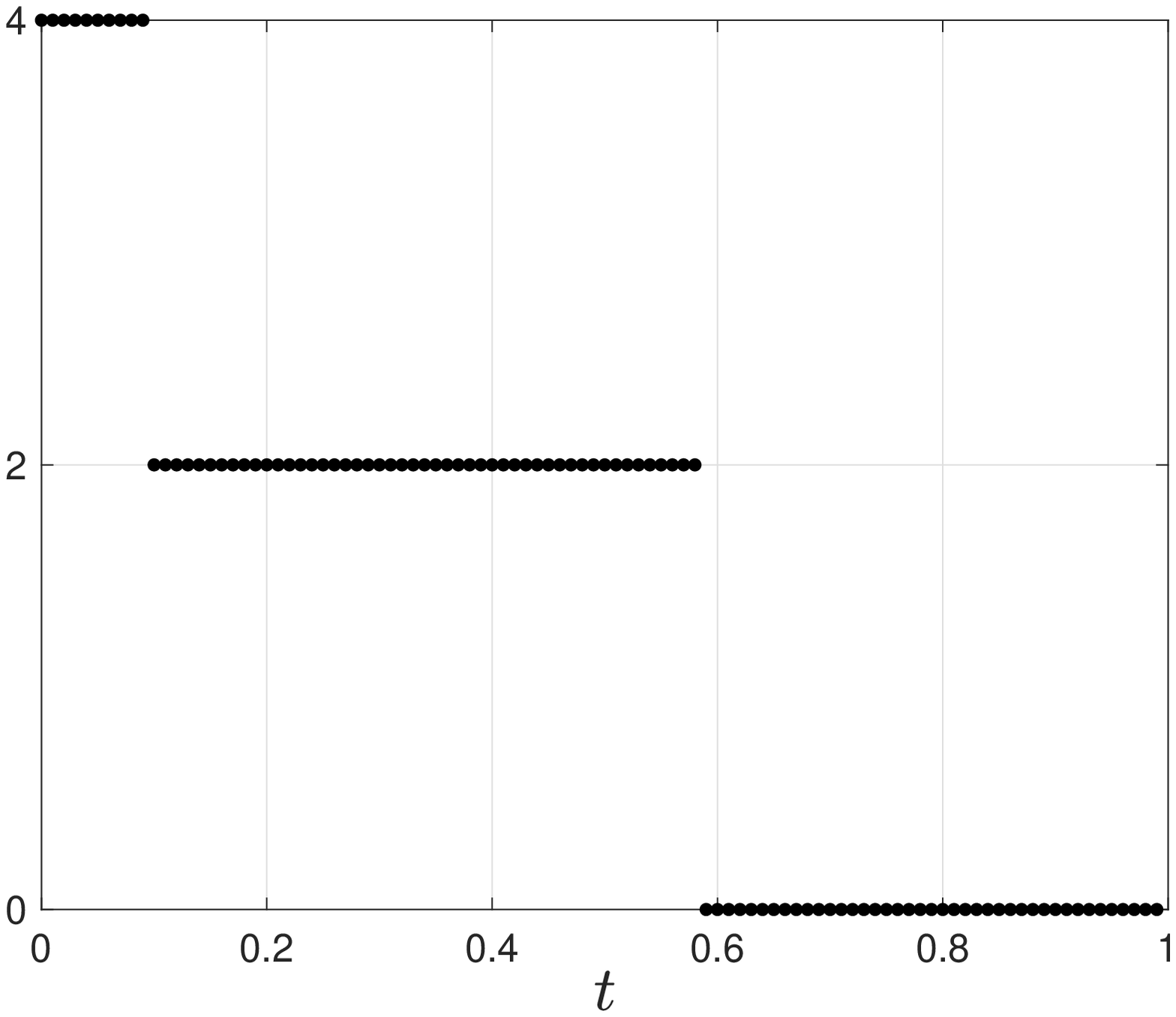}
\end{center}
\end{subfigure}
\hfill
\begin{subfigure}[t]{0.5\textwidth}
\begin{center}
\includegraphics[width= 8cm,height=7cm]{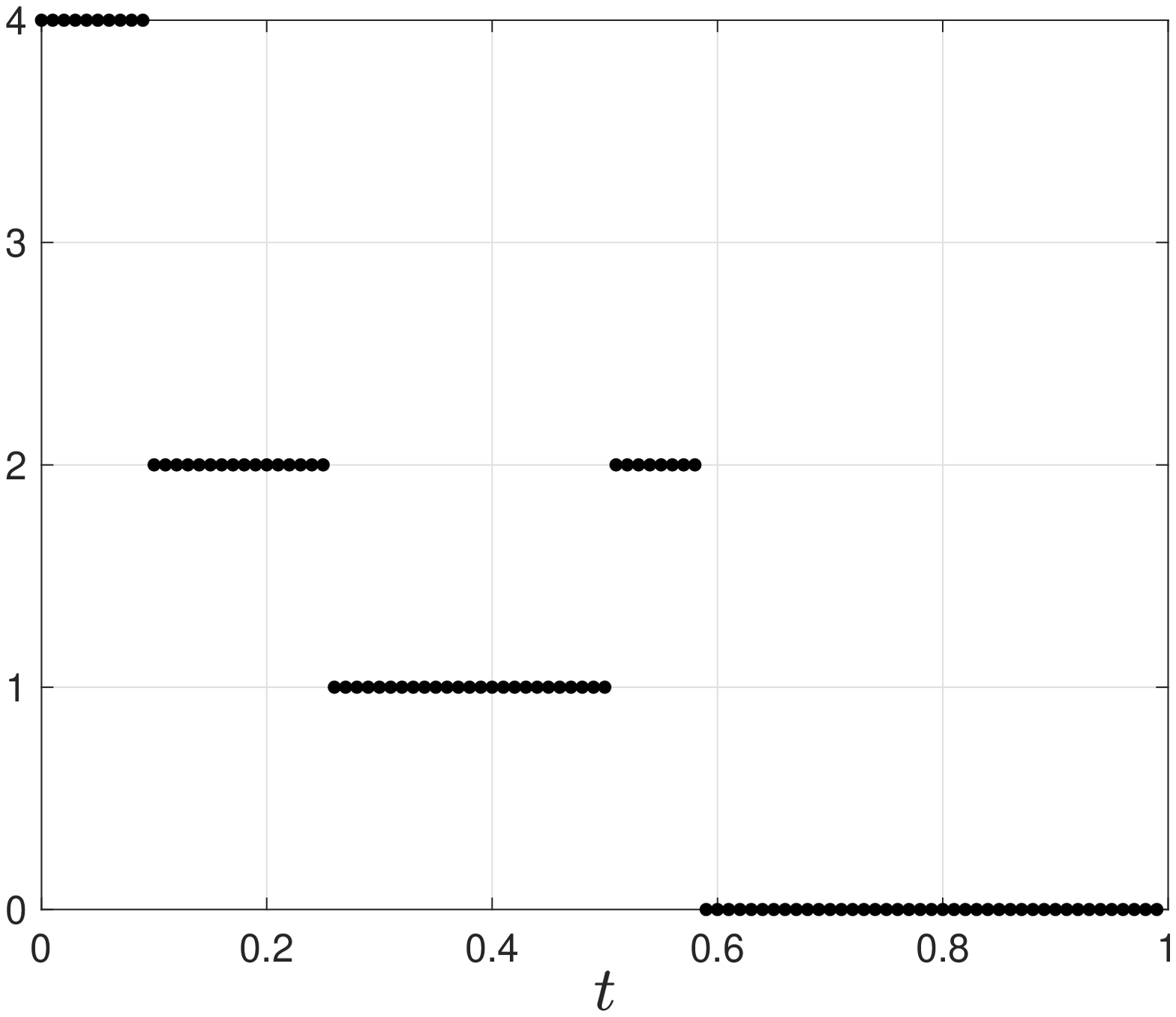}
\end{center}
\end{subfigure}
\caption{$s^{-}_c(x(t))$ (left) and~$s^{-} (x(t)) $ (right) 
as a function of $t\in[0,1]$ for the system in Example~\ref{exp:lin_n5}.\label{fig:sands}}
\end{figure*}

\begin{Example}\label{exp:lin_n5}
Consider the system $\dot x = Ax$, with
\[
A=\begin{bmatrix}
-4 &   1 &         0   &      0  &  0 \\
    2 &   -4 &    4 &         0  &       0\\
         0 &    7/2 &   -4 &    5/2 &         0\\
         0    &     0  &  0 &   -4 &    1 \\
    5/4 &        0       &  0 &   3/2 &   -4 \\
\end{bmatrix}.
\]
Note that $A\in (\Q^+ \setminus \M^+)$, thus the system is {\model} but not TPDS.
 Fig.~\ref{fig:sands} depicts $s_c^{-}(x(t))$ and~$s^{-}(x(t))$ as a function of $t\in[0,1]$ for the initial condition
$x(0)=\begin{bmatrix} -0.6407  &  1.8089   & -1.0799   &  0.1992  &  -1.5210 \end{bmatrix}'$. It may be noticed that $s_c^{-}(x(t))$ is piecewise-constant and that at any point where its value changes it decreases by two.
On the other-hand 
 $s^{-}(x(t))$ both decreases and increases (as the system is not~TPDS), and the
 increase agrees with the structure described in Remark~\ref{rem:QnotM}.
\end{Example}



The proof of Thm..~\ref{thm:mecsauc} is based on analyzing the dynamics of all 
the odd minors of~$\Phi(t,t_0)$. To do that, we briefly review multiplicative and additive compound matrices (see, e.g.,~\cite{muldo1990,fulltppaper}).
Given~$A\in\R^{n\times n}$ and~$p\in\{1,\dots,n\}$, consider the
$\binom{n}{p}^2$
 minors of~$A$ of size~$p\times p$. 
Recall that each such  minor is defined by a set of row indexes~$1\leq i_1<i_2<\dots<i_p\leq n$ and column indexes~$1\leq j_1<j_2<\dots<j_p\leq n$, and is denoted by~$A(\alpha|\beta)$, where~$\alpha:=\{i_1,\dots,i_p\}$ and~$\beta:=\{j_1,\dots,j_p\}$.
The~$p$th \emph{multiplicative  compound matrix}~$A^{(p)}$ is the~$\binom{n}{p}\times  \binom{n}{p}$ matrix
that 
includes all these minors ordered lexicographically.
The Cauchy-Binet formula yields
\be\label{eq:miucpxdw}
(AB)^{(p)}=A^{(p)}B^{(p)},
\ee
 justifying the term multiplicative compound.

The $p$th \emph{additive compound matrix} of~$A$
is   defined  by
\be\label{eq:defapp}
				A^{[p]}:= \frac{d}{dh}  (I+hA)^{(p)} |_{h=0}.
\ee
In other words,~$A^{[p]}$ is the the term that multiplies~$h$ 
in the Taylor series expansion of~$(I+hA)^{(p)}  $.  
Using~\eqref{eq:defapp} and~\eqref{eq:miucpxdw} gives~$(A+B)^{[p]}=A^{[p]}+B^{[p]}$, justifying
 the term additive compound.

It can be shown~\cite{schwarz1970}
 that if~$\Phi$ satisfies~\eqref{eq:mdfe} then for any~$p\in\{1,\dots,n\}$,
\be\label{eq:povt}
\frac{d}{dt} \Phi^{(p)} =A^{[p]} \Phi^{(p)}.
\ee
Thus, the dynamics of~$\Phi^{(p)} $, i.e. the dynamics of the minors of order~$p$
of~$\Phi $, is also a linear system  with the matrix~$A^{[p]}$.

The matrix~$A^{[p]}$ can be determined explicitly. 
The entry of~$A^{[p]}$ corresponding to~$(\alpha|\beta)=(i_1,\dots,i_p|j_1,\dots,j_p) $  is:
\be\label{eq:poltr}
	\begin{cases}
			 \sum_{k=1}^p a_{i_k i_k},		&\text{if }   |\alpha\cap \beta|=p,\\
					(-1)^{\ell+m} a_{i_\ell j_m},		& \text{if } |\alpha\cap \beta|=p-1  \text{ and } i_\ell\not = j_m,\\
					0,& \text{otherwise.} 
\end{cases}
\ee
The first line  in~\eqref{eq:poltr} corresponds to the case 
where~$i_k=j_k$ for all~$k=1,\dots,p$, i.e. to the diagonal entries of~$A^{[p]}$.  
The second line  describes the case where   all the indexes in~$\alpha$ and~$\beta$ coincide, except for a single index~$i_\ell\not = j_m $. 
Eq.~\eqref{eq:poltr} is  usually proven by manipulating determinants~\cite{schwarz1970}
or using exterior powers~\cite{fiedler_book}.

For example, consider the case~$A=\{a_{ij}\}_{i,j=1}^4$. 
Then~\eqref{eq:poltr} yields
$
A^{[1]}=A,
$
  \\
\be\label{eq:a2445}
A^{[2]}=\begin{bmatrix} 
																a_{11}+a_{22} & a_{23} & a_{24} & -a_{13} & -a_{14} & 0 \\
																a_{32}      &a_{11}+a_{33} & a_{34} & a_{12} & 0 & -a_{14} \\   
																a_{42}     &a_{43}   &a_{11}+a_{44} & 0 & a_{12} & a_{13} \\   
																-a_{31}      &a_{21} & 0 & a_{22}+a_{33} & a_{34} &  -a_{24} \\
																-a_{41}      &0 & a_{21} & a_{43} & a_{22}+a_{44} &  a_{23} \\ 
																0      &-a_{41} & a_{31} & -a_{42} & a_{32} &  a_{33}+a_{44}  
\end{bmatrix},
\ee
\be\label{eq:a3445}
A^{[3]}=\begin{bmatrix}     
																a_{11}+a_{22}+a_{33} & a_{34} & -a_{24} & a_{14} \\
																a_{43}& a_{11}+a_{22}+a_{44} & a_{23} & -a_{13}  \\
																-a_{42}& a_{32} & a_{11}+a_{33}+a_{44} & a_{12}  \\
                                a_{41}& -a_{31}& a_{21}& a_{22}+a_{33}+a_{44} 
 \end{bmatrix},
\ee \\
and~$A^{[4]}=\tr(A)$. 

 
{\sl Proof of Thm.~\ref{thm:mecsauc}.}
  Suppose that
	\be\label{eq:assumpq}
	A\in \Q^+.
	\ee
	We will show that for any odd~$p$ the matrix~$A^{[p]}$ is Metzler and irreducible.
	For~$p=1$ this is immediate, as~$A^{[1]}=A$.
	
	Pick an odd~$p\in\{3,\dots,n\}$. By~\eqref{eq:poltr},
	every off-diagonal entry of~$A^{[p]}$ is either zero  
	or has the form~$(-1)^{\ell+m} a_{i_\ell j_m}$. 
	Since~$A\in\Q^+$, $A ^{[p]}$ can have a non-zero
	off-diagonal   entry  only 
	in the following three cases.

\noindent {\sl Case 1.} The entry~$a_{i_\ell j_m}$ is
on the super- or sub-diagonal of~$A$.
Then~$\{i_\ell,j_m\}=\{k,k+1\}$ for some~$k\in\{1,\dots,n-1\}$. This means that $\alpha$ and~$\beta$ 
share $p-1$ common indexes, and one of $\{k,k+1\}$ is in $\alpha$ (and not in $\beta$) and the other in $\beta$ (and not in $\alpha$). 
This implies that $\ell=m$, and thus $(-1)^{\ell+m} a_{i_\ell j_m}=(-1)^{2\ell}a_{i_\ell j_\ell}$. Since~$i_\ell \not = j_\ell$ and~$A$ is Metzler,
we conclude that~$(-1)^{\ell+m} a_{i_\ell j_m}\geq 0$.
 
\noindent {\sl Case 2.}  
The entry~$a_{i_\ell j_m}$ is
 entry $(1,n)$ in $A$. 
This means that~$\ell=1$, $m=p$, so
$(-1)^{\ell+m} a_{i_\ell j_m}=(-1)^{p+1}a_{1n}\geq 0$, as~$p$ is odd. 

\noindent {\sl Case 3.} 
The entry~$a_{i_\ell j_m}$ is
 entry $(n,1)$ in $A$. 
 Then~$\ell=p$, $m=1$, so~$(-1)^{\ell+m} a_{i_\ell j_m}=(-1)^{1+p}a_{n1}\geq 0$.

We conclude that~$A^{[p]}$ is Metzler.
We now show that~\eqref{eq:assumpq} implies that~$A^{[p]}$ is irreducible.
We   introduce some notation. 
Let~$G^{[p]}$ denote the  adjacency  graph  associated with
the matrix~$A^{[p]}$. Every node in this graph corresponds to
 a set of~$p$    indexes~$1\leq i_1<i_2<\dots<i_p\leq n$, and there are~$\binom{n}{p}$ nodes. 
There is a  directed edge 
from node~$\alpha=\{i_1,\dots,i_p\}$ to~$\beta=\{j_1,\dots,j_p\}$
  if exactly~$p-1$ entries of~$\alpha$ and~$\beta$   coincide,    
the  two remaining indexes 
 are~$i_\ell\not =j_m$, and~$a_{i_\ell j_m}\not =0$.  
We write~$\alpha \leq \beta$ if~$i_k\leq j_k$ for all~$k\in \{ 1,\dots,p\}$.  
We write~$\alpha \leadsto \beta$ if there is a path in~$G^{[p]}$  
starting at~$\alpha$ and ending at~$\beta$. 

We consider two cases.

\noindent {\sl Case 1.} Suppose that all the entries on the sub- and super-diagonals of~$A$ are positive. 
  It follows from the results of
Schwarz~\cite{schwarz1970}, that considered the case~$A\in \M^+$,   that in this case~$G^{[k]}$
is strongly connected for all~$k$, so~$A^{[k]}$ is irreducible for all~$k$. In particular,~$A^{[p]}$ is irreducible. 

\noindent {\sl Case 2.} Suppose that~$A\in \Q^+$, but  
there exist~$i,j$, with~$|i-j|=1$ such that~$a_{ij}=0$. 
We may assume that for some~$i\in\{2,\dots,n\}$ entry~$(i,i-1)$ of~$A$ is zero.
We know that~$A$ is irreducible, and using the structure of~$\Q^+$ 
this implies that
\be\label{eq:aaiii}
a_{12}>0, \; a_{23}>0,\dots,a_{n-1n}>0,\; a_{n1}>0 .
\ee 
Let
\begin{align*}
\underline \gamma&:=\{1,2,\dots,p\},\\
\bar \gamma & : = \{ n-p+1  , n-p+2, \dots,  n-1,n \} ,
\end{align*}
 i.e. the lowest and highest nodes in
 the lexicographic ordering, respectively.  
It follows from~\eqref{eq:poltr} and~\eqref{eq:aaiii} that for any two
  nodes~$\alpha \not = \beta$ satisfying~$\alpha\leq \beta$, we have~$\alpha \leadsto \beta$.
In particular, there is a path from~$\underline \gamma$ to any node, and there is a path from any node to~$\bar \gamma$. 

Since~$a_{n1}>0$, 
the   edges 
\begin{align*}
 \{i_1,i_2,\dots,i_{p-1} ,n\} &\to  \{1,i_1,i_2,\dots,i_{p-1}  \} , 
\end{align*}
with~$2\leq i_1<i_2<\dots<i_{p-1}<n$ are included in~$G^{[p]}$.
In particular, this includes the edges
$
 \bar \gamma   \to  \delta:= \{1, n-p+1  , n-p+2,
\dots, n-2, n-1  \} $, and~$\zeta:=\{2,3,\dots,p ,n\}  \to \underline  \gamma$.

Now pick   two distinct nodes~$\alpha$, $\beta$. 
We will show that~$ \alpha \leadsto \beta $.
If
\[
p\geq n-2
\]
then~$\delta 
				\leq \zeta$, so~$G^{[p]}$ includes the path:
\[
				\alpha \leadsto \bar \gamma \to \delta  
				\leadsto \zeta \to \underline \gamma \leadsto \beta.
\]
If
\[
p\geq n-3
\]
then~$\{ 1,2,n-p+1,n-p+2,\dots, n-2   \} 
				\leq \zeta$, so~$G^{[p]}$ includes the path:
\begin{align*}
				\alpha &\leadsto \bar \gamma \to \delta
				\leadsto \{ 2,n-p+1,n-p+2,\dots, n-2,n  \} \\ 
				&\to \{ 1,2,n-p+1,n-p+2,\dots, n-2   \} \leadsto \zeta
				\leadsto \underline \gamma \leadsto \beta.
\end{align*}
Proceeding in this fashion we see that there is always
 a path from~$\alpha$ to~$\beta$. Thus,~$G^{[p]}$ is strongly connected, so~$A^{[p]}$ is irreducible.

Summarizing, for any odd~$p$ we have that~$A\in\Q^+$ implies that~$A^{[p]}$ is Metzler and irreducible. This means that~\eqref{eq:povt} is a strongly cooperative dynamical system~\cite{hlsmith}, and since~$\Phi^{(p)}(t_0,t_0)=I$,
 we conclude that~$\Phi^{(p)}(t,t_0)$ is 
componentwise positive for all~$t>t_0$. In other words, 
all minors of order~$p$
are positive for all~$t>t_0$, so the system is~{\model}.

To prove the converse implication, assume that~$A\not \in \Q^+$.
Then one of the following three cases holds.

\noindent {\sl Case 1.} 
The matrix~$A$ is not Metzler. 
Then~$a_{ij}<0$ for some~$i\not =j$.
 Since~$\Phi(t_0)=I$ 
and
\be\label{eq:alfig}
\dot \Phi(t_0)=A \Phi(t_0)=A, 
\ee
it follows that   entry~$(i,j)$ of~$\Phi (t_0+\varepsilon)$ is negative for 
all~$\varepsilon>0$ sufficiently small. Thus, the system is not~{\model}.

\noindent {\sl Case 2.} 
The matrix~$A$ is Metzler, but not irreducible. In this case, it  is straightforward to show 
that there exists an entry of~$\Phi(t)$ that is zero for all~$d:=t-t_0>0$ sufficiently small. 
 Thus, the system is not~{\model}.

\noindent {\sl Case 3.} 
The matrix~$A$ is Metzler,  irreducible, but 
there exist indices~$w,q\in\{1,\dots,n\}$ such that
 $a_{wq}\ne 0$ with 
\be\label{eq:wqe}
1< |w-q|  \text{ and }  |w-q|<n-1.
\ee
The first inequality here means that~$a_{wq}$ is not on the main, super- or sub-diagonal.
The second inequality implies that~$a_{wq}$ is not   entry~$(1,n)$ nor~$(n,1)$.
We will show that   there exists an odd~$p$ 
 such~$A^{[P]}$ is not Metzler.
 If~$a_{wq}<0$ then~$A^{[1]}=A$ is not Metzler. 
Thus, it is enough to
 consider the case~$a_{wq}>0$.
 Assume, without loss of generality, that $w<q$ (the analysis in the case~$w>q$ is similar). It follows from~\eqref{eq:wqe}
that  there exists~$k\in\{2,\dots,n-1\}$ such that $w<k<q$. We consider two cases.

\noindent {\sl Case 1.} Suppose that~$q\ne n$. Let $\alpha:=\{w,k,n\}$ and $\beta:=\{k,q,n\}$. Then $|\alpha \cap \beta|=2$, and $i_1=w\ne j_2=q$. Thus, the entry in~$A^{[3]}$  corresponding to~$(\alpha,\beta)$ is $(-1)^{\ell+m} a_{i_\ell j_m}=(-1)^{1+2}a_{wq}<0$, so $A^{[3]}$ is not Metzler.

\noindent {\sl Case 2.} Suppose that~$q=n$. Note that since $q-w<n-1$, $w>1$. Let~$\alpha:=\{1,w,k\}$ and $\beta:=\{1,k,q\}$. Then $|\alpha \cap \beta|=2$, and $i_2=w\ne j_3=q=n$. Thus, the entry in~$A^{[3]}$ corresponding  to~$(\alpha,\beta)$  is $(-1)^{\ell+m} a_{i_\ell j_m}=(-1)^{2+3}a_{wn}<0$, so again $A^{[3]}$ is not Metzler.

We conclude that if $A\notin \Q^+$ then at least one 
of the matrices~$A^{[1]}$, $A^{[3]}$ is either not Metzler or reducible.
Since~$\Phi(t_0 )=I$,~$\Phi^{(p)}( t_0)=I$ for all~$p$. 
If~$A^{[p]}$ is not Metzler then there exist~$i,j$, with~$i\not =j$,
such that entry~$(i,j)$ of~$A^{[p]}$ is negative.
Since~$\dot {\Phi}^{(p)}(t_0)=A^{[p]} \Phi^{(p)}(t_0)=A^{[p]}$, 
it follows that entry~$(i,j)$ of~$\dot \Phi^{(p)}(t_0)$ is negative 
 and thus entry~$(i,j)$ of~$\Phi^{(p)}(t_0+\varepsilon)$ is negative for all~$\varepsilon>0$ sufficiently small. 
Thus, the system is not~{\model}.
 This completes the proof of 
 Thm.~\ref{thm:mecsauc}.~\hfill{$\square$} 

\begin{Remark}
Thm.~\ref{thm:mecsauc} 
implies in particular that if~$A\in M^+$  then~$\dot x=Ax$ is~{\model}. 
In fact in this case  the system is also~TPDS, so its transition matrix satisfies the~SVDP and thus it satisfies the~SCVDP. 
In this respect,~TPDS is a special case of~{\model}.
\end{Remark}

\begin{Remark}\label{rmk:A1A3}
The proof of Thm.~\ref{thm:mecsauc} shows that if~$A\in\Q^+$
then the system is {\model} and if~$A\not\in\Q^+$ 
then at least one of~$A^{[1]}=A$, $A^{[3]}$ is not Metzler or irreducible. 
Thus, for any dimension~$n$ it is sufficient to verify that~$A$ 
and~$A^{[3]}$ are Metzler and irreducible in order to establish {\model}. 
\end{Remark}

\begin{figure*}
\begin{center}
\includegraphics[width= 12cm,height=4cm]{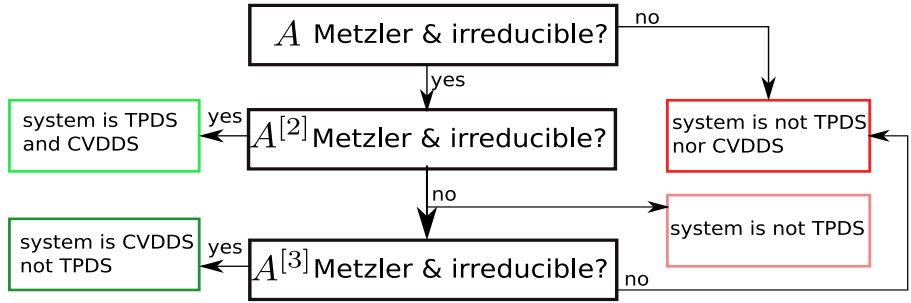}
\caption{Flowchart for checking if~$\dot x(t) = Ax(t)$ is {\model} 
and/or~TPDS.}\label{fig:is_ctpds_tpds}
\end{center}
\end{figure*}

Combining Remark~\ref{rmk:A1A3} with the results in~\cite{fulltppaper} yields
 a simple flowchart for establishing if a system~$\dot x(t) = Ax(t)$ is {\model} and/or TPDS.
 This is depicted in Fig.~\ref{fig:is_ctpds_tpds}. 

\begin{Example}
Consider the case~$A=\{a_{ij}\}_{i,j=1}^4$. 
In this case  the matrices~$A^{[2]}$, $A^{[3]}$ are given in~\eqref{eq:a2445} 
and~\eqref{eq:a3445}
and~$A^{[4]}=\tr(A)$.
 To guarantee that~$A^{[1]}=A$ is Metzler,
$a_{ij}$ must be nonnegative for all~$i\not =j$. 
Now~\eqref{eq:a3445} shows that 
to guarantee that~$A^{[3]}$ is Metzler,
we must have~$a_{13}=a_{31}=a_{24}=a_{42}=0$, 
so~$A \in \Q^+$. 
\end{Example}

if~$\dot x=Ax$    is~{\model} then we know from  Thm.~\ref{thm:exp_is_ctp}
    that for any initial condition~$x(0)\not = 0 $
there exists a time~$T\geq 0$ such that~$s^-_c(x(t))=s^+_c(x(t))$ for all~$t\geq T$. 
However, it is important to note that this does 
not necessarily imply that after some finite time
  the state-variables  do not change their
	sign anymore. The next example demonstrates this. 
\begin{Example}\label{exa:unbou}
Consider~$\dot x=Ax$ with~$n=3$ and~$A=\begin{bmatrix} 0&0&1\\ 1&0&0\\ 0&1&0\end{bmatrix}$. 
This is {\model} as~$A\in\Q^+$. The solution for~$x(0)=\begin{bmatrix} 1&-2&1\end{bmatrix}'$ is 
\[
x(t)=\exp(-t/2)
\begin{bmatrix}  
				\cos(pt/2)+p\sin(pt/2) \\ -2\cos(pt/2) \\ \cos(pt/2)-p\sin(pt/2)													
\end{bmatrix},
\]
where~$p:=\sqrt{3}$. This means that \emph{every}~$x_i(t)$ changes sign an unbounded 
 number of times. Note that here~$s^-_c(x(t))=s^+_c(x(t))=2$
 for all~$t\geq 0$.
\end{Example}

We now turn to
consider~\eqref{eq:ltv}  (and~\eqref{eq:mdfe}) with~$A(t)$   time-varying.  
\subsection{The case~$A(t)$ continuous}
Suppose that~$t\to A(t)$ is a continuous matrix function
on a time interval~$(a,b)$.
To derive a necessary and sufficient condition for~{\model}, we 
 first introduce more notation.
For a matrix~$Q$ we write~$Q \geq 0$ [$Q\gg 0$] if every entry of~$Q$ is nonnegative
[positive]. The next result provides a necessary and
sufficient condition for the solution 
 of a linear matrix differential equation
 to be componentwise  positive for all time.
\begin{Proposition}\label{prop:resp}
Let~$A(\cdot):(a,b)\to\R^{n\times n}$ be  a continuous matrix function,
and assume that~$A(\tau)$   is Metzler for all~$\tau\in(a,b)$.
Then the following two conditions are equivalent.
\begin{enumerate}[(1)] 
\item \label{cond:vbd} for  any interval~$[p,q]$, with~$a<p<q<b$, there exists~$t^* \in [p,q]$ such that~$A(t^*)$ is   irreducible;
\item \label{cond:zety} for any pair~$(t_0,t)$ with~$a<t_0<t<b$ the solution 
 of~$\dot \Phi(\tau)=A(\tau)\Phi(\tau)$, $\Phi(t_0)=I$, 
satisfies~$\Phi(t)   \gg 0$.  
\end{enumerate} 
\end{Proposition}

{\sl Proof of Proposition~\ref{prop:resp}.}
Throughout we will use  the well-known graph-theoretic interpretation  of irreducibility, namely, 
that   a matrix~$B\in \R^{n\times n}_+$ (and thus also a Metzler matrix~$B$)  is irreducible iff
 its associated adjacency graph is strictly connected~\cite[Chapter~6]{matrx_ana}.
 This  immediately implies that if~$A(\tau)$ is 
  irreducible   and~$A(\cdot)$ is continuous  then there exists a time interval containing~$\tau$ such that~$A(s)$
 is irreducible for all~$s$ in this time interval. 

For~$s\in(a,b)$, 
let~$E_+(s)$ denote the edges in the adjacency matrix with positive weights at time~$s$ (corresponding to 
positive off-diagonal  entries of~$A(s)$), and 
let~$E_0(s)$ denote the edges in the adjacency matrix with zero weights at time~$s$ (corresponding to 
zero off-diagonal  entries of~$A(s)$).

Suppose      that condition~\eqref{cond:vbd} in Proposition~\ref{prop:resp} holds. 
Fix~$t_0,t$ with~$a<t_0<t<b$. 
By condition~\eqref{cond:vbd}, there exists~$t^* \in [t_0,t]$ such that~$A(t^*)$ is irreducible. 
By continuity, there exists a time interval~$L \subseteq [t_0,t]$ such that~$A(s)$ is irreducible for all~$s\in L$ and, furthermore, the matrix
\[
\underline A:=\min_{s\in L }A(s)
\]
 is Metzler and irreducible. For any~$\tau \in (a,b)$, let
\[
B(\tau):=\begin{cases}
											\underline A,& \tau\in L,\\
											A(\tau),& \text{otherwise},
\end{cases}
\]
and let~$\Omega$ denote the solution of~$\dot \Omega=B \Omega $, with~$\Omega(t_0)=I$.
By known results on strongly cooperative dynamical systems (see, e.g.~\cite{fulltppaper}),
~$\Omega(t) \gg 0$. Let~$\Delta:=\Phi-  \Omega$. Then 
 \[
				\dot \Delta = B \Delta  +  (A-B) \Phi,
\]
and this
implies that~$\Delta(t)\geq 0$, so~$\Phi(t)\geq \Omega(t)\gg 0$.

To prove the converse implication, assume  that condition~\eqref{cond:vbd} does \emph{not} hold. 
Thus,  there exist~$p,q$, with~$a<p<q<b$, such that~$A(\tau)$ is    
reducible for \emph{all}~$\tau\in[p,q]$. 
By continuity, there exists~$\varepsilon>0$ sufficiently small such that the
 edges in~$E_+(p)$ belong to~$E_+(s)$ for all~$s\in [p,p+\varepsilon]$. 
Since~$A(\tau)$ is    
reducible for all~$\tau\in[p,p+\varepsilon]$, there exists a non empty set of
 edges that are zero for all~$\tau\in[p,p+\varepsilon]$ and when these edges are zero the associated graph is not strictly connected. 
Hence,  there exists a permutation matrix~$P\in\{0,1\}^{n\times n}$ such that   
\[
PA(s)P'=\begin{bmatrix}   B(s)& C(s) \\ 0& D(s)  \end{bmatrix},\quad \text{ for all } s\in[p,p+\varepsilon],
\]
where the~$0$ denotes an~$(n-r)\times r$ zero matrix, with~$ 1\leq r \leq n-1$.
This implies that condition~\eqref{cond:zety} does not hold for~$t_0=p$ 
and any~$t
\in (p,p+\varepsilon]$.
 This completes the proof.~\hfill{$\square$}

We can now state
 the main result in this subsection.
Let~$\Q \subset \R^{n\times n}$ denote the set of~$n\times n$
matrices that   are Metzler  and can have    non-zero entries   only
  on the main diagonal, super- and sub-diagonals and entries $(1,n)$
and~$(n,1)$. For example, for~$n=4$ any  matrix in the form:
\[
\begin{bmatrix}
									*& +& 0&0 \\ 0& *& 0&0 \\ 0& 0& *&+ \\ + & 0& 0&*
\end{bmatrix},
\]
where~$*$ denotes some value and~$+$ denotes a positive value, 
is in~$\Q $. 
\begin{Theorem}\label{thm:nsct}
					Suppose that~$A(\cdot)$ is a continuous matrix function on~$(a,b)$. 
					Then the following conditions are equivalent:
					\begin{enumerate}[(1)]
					\item \label{enu:poiid} $A(s) \in \Q$ for all~$s\in (a,b)$, and 
					for  any interval~$[p,q]$, with~$a<p<q<b$, there exists~$t^* \in [p,q]$ such that~$A(t^*)$ is   irreducible;
					\item  \label{enu:eccsd} the system~\eqref{eq:mdfe} is {\model}.
           \end{enumerate} 
\end{Theorem}
Note that this generalizes Theorem~\ref{thm:mecsauc}. Indeed, if~$A$ is a
 constant matrix then condition~\eqref{enu:poiid}
means that~$A\in   Q$ and~$A$ is irreducible, so~$A\in \Q^+$. 

{\sl Proof of Theorem~\ref{thm:nsct}.}
If~$A(s) \not \in \Q$ for some~$s \in (a,b)$ then using
 continuity implies that~$A(\cdot)  \not \in \Q$ on a time interval that includes~$s$
and arguing as in the proof of Theorem~\ref{thm:mecsauc} implies
 that~\eqref{eq:mdfe} is not~{\model}. If~$A(\cdot)$ is reducible on some
 time interval~$[p,q]$ 
then it follows from Proposition~\ref{prop:resp} 
that~\eqref{eq:mdfe} is not~{\model}.

To prove the converse implication, assume that condition~\eqref{enu:poiid} holds.
Pick~$(t_0,t)$ such that~$a<t_0<t<b$, and an odd integer~$p$. 
Arguing as in the proof of   Theorem~\ref{thm:mecsauc} shows that~$A^{[p]}$ 
is also  Metzler and  furthermore that
  there exists~$t^* \in [t_0,t]$ such that~$A^{[p]}(t^*)$ is irreducible. Using Proposition~\ref{prop:resp}  
implies that the solution of~$\dot \Phi^{(p)}= A^{[p]} \Phi^{(p)}$, with~$\Phi^{(p)}(t_0)=I$,
 satisfies~$\Phi^{(p)}(t) \gg 0$, i.e. all minors of order~$p$ are positive at time~$t$. 
Since this holds for any odd~$p$, this completes the proof.~\hfill{$\square$}

\subsection{The case~$A(t)$ measurable}
In this subsection,  we assume that
\begin{align}\label{eq:sas}
\text{ }  A:(a,b) \to \R^{n\times n} 
&\text{ is a  matrix of locally (essentially)}\nonumber \\ & \text{ bounded  measurable 
functions.}
\end{align}
This   general case  is important in the context of control systems. Indeed,  
 consider~$\dot x=f(x,u)$, with~$u$ a control input. The associated
 variational equation
 is~$\dot z=J(x,u)z$, where~$J=\frac{\partial f}{\partial x}$ is the Jacobian of~$f$. In many cases,
 for example when considering optimal controls, one must allow measurable controls (see, e.g.~\cite{liberzon_opt_cont})  
and thus~$t\to J(t)$ is typically a measurable, but  
 not necessarily continuous, matrix function. 

 It is well-known 
 that~\eqref{eq:sas} implies 
that~\eqref{eq:mdfe} admits a unique, locally absolutely continuous,    nonsingular 
solution for   all~$t \in(a,b)$ 
(see, e.g.,~\cite[Appendix~C]{sontag_book}). 
The next result provides a sufficient condition for~\eqref{eq:ltv} to be~{\model}.

 \begin{Theorem} \label{thm:meas}
Suppose that~\eqref{eq:sas} holds with~$A(t)\in \Q^+$ for almost all $t\in(a,b)$. Furthermore, suppose that there exists~$\delta>0$
 such that all  the non-zero  entries of~$A(t)$ are larger or equal to $\delta$ for almost all $t$. 
 Then~\eqref{eq:mdfe} is~{\model}.
\end{Theorem}

{\sl Proof of Thm.~\ref{thm:meas}.}
Pick an odd $p\in\{1,\dots,n\}$. Consider the dynamics~$\dot \Phi^{(p)}(t)= A^{[p]}(t) \Phi^{(p)}(t)$. 
We know from the analysis in the proof of Thm.~\ref{thm:mecsauc} that
every off diagonal entry of $A^{[p]}(t)$ is either zero or equal or larger than $\delta$ for almost all $t\in(a,b)$. 
To analyze~$\Phi^{(p)}(t)$, pick an arbitrary~$1\leq k \leq \binom{n}{p}$, and let~$u(t)$
 denote the~$k$th column of~$\Phi^{(p)}(t)$. Then~$\dot u(t)=A^{[p]}
(t)u(t)$,
with~$u(t_0)=e^k$, where~$e^k$ is the~$k$th canonical vector in~$\R^{\binom{n}{p}}$.
For any~$j\in\{1,\dots,\binom{n}{p}\}$ the linear equation for~$\dot u_j(t)$
implies that there exists~$c_j \in \R$ such that~$u_j(t)\geq \exp(c_j (t-t_0))u_j(t_0)$ for all~$t\geq t_0$.
 In particular, if~$u_j(\tau)>0$ at some time~$\tau$ then~$u_j(t)>0$ for all~$t\geq \tau$. Thus, $u_k(t)>0$ for all~$t\geq t_0$.
Pick a time~$\tau \geq t_0$, and let~$w\geq 1$ denote the number of entries~$j$
such that~$u_j(\tau)>0$. Without loss of generality, assume 
that these entries are~$j=1,2,\dots, w$. Write~$\dot u=A^{[p]} u$   as 
\[
			\dot u= \begin{bmatrix} *&*\\G& *\end{bmatrix} \begin{bmatrix} u_1\\ \vdots \\
			u_w\\0    \end{bmatrix},
\]
where the~$0$ denotes a vector of~$s:=\binom{n}{p}-w$ zeros, and~$G(t)\in \R^{s \times w} $. 
Since~$A^{[p]}(t)$ is irreducible for almost all~$t$, $G(t)$ has a non-zero
 entry that is larger or equal to~$\delta$ for almost all~$t$. 
This implies that at least~$w+1$ entries in $u_j(t)$ are positive for all~$t>\tau$.
Our assumption on~$A(t)$ implies that we can now use an inductive argument
to conclude that all the  entries of~$u(t)$ are positive for all~$t>t_0$. Since this 
holds for any odd~$p$   and any index~$k$,
we conclude that every odd minor  of~$\Phi(t)$ is positive for all~$t>t_0$.
 This completes  the proof of Thm.~\ref{thm:meas}.~\hfill{$\square$}

In the remainder of this section, we demonstrate  how 
  Theorem~\ref{thm:nsct} can be used to show that the variational equation associated with a 
	nonlinear dynamical model, called the  \emph{ribosome  flow model on a ring}~(RFMR),
	is~{\model}.

There is considerable recent interest 
in understanding the dynamics of mRNA translation. This is due to:
 (1)~the introduction of several novel    methods allowing
 to track  mRNA translation  in the cell~\cite{natrevtrans,INGOLIA201622}; and  
(2)~a growing understanding that regulation of mRNA translation 
plays an important role in the 
precise tuning of the expression of each gene in the genome, and
that this is
 critical for many aspects of cell function~\cite{pnas_sh,Millseaan2755,slow_rib}.

The~RFMR has been used to study the flow of ribosomes along a circular mRNA during the process of 
translation~\cite{RFMR,alexander2017,RFM_r_max_density}. 
This model    includes~$n$ consecutive sites located along a circular ring.
 The normalized occupancy level (or density) of site~$i$ at time~$t$
    is described
		by a state variable~$x_i(t): \R_+ \to [0,1]$, $i=1,\dots,n$,  where~$x_i(t)=0$ [$x_i(t)=1$] means that site~$i$ is completely free [full].  
The transition between sites~$i$ and site~$i+1$ 
is regulated by a parameter~$\lambda_i>0$. 
 
From here on we interpret  all indexes modulo~$n$, i.e.~$x_0=x_n$ and~$x_{n+1}=x_1$.
The dynamics of the~RFMR   is given by $n$ nonlinear first-order ordinary differential equations:
\be\label{eq:rfmr_all}
\dot{x}_i=\lambda_{i-1}x_{i-1}(1-x_i)-\lambda_i x_i(1-x_{i+1}),\quad i=1,\dots,n.
\ee
 This  can be explained as follows. The flow of particles from site~$i$ to site~$i+1$ is~$\lambda_{i} x_{i}(t)
(1 - x_{i+1}(t) )$. This flow is proportional to $x_i(t)$, i.e. it increases
with the occupancy level at site~$i$, and to $(1-x_{i+1}(t))$, i.e. it decreases as site~$i+1$ becomes fuller.  
Note that the maximal possible  flow  from site~$i$ to site~$i+1$  is the transition rate~$\lambda_i$.
Eq.~\eqref{eq:rfmr_all} thus states that the change in the  state-variable~$x_i$ as a function of 
time equals the  flow entering site~$i$ from site~$i-1$,    
 minus the flow exiting site~$i$ to site~$i+1$. 

It is not difficult to show that~$[0,1]^n$ is an invariant set for the~RFMR. 
It is easy to verify that~$1_n$ and~$0_n$ are equilibrium points of the~RFMR
(corresponding to the case where \emph{all} sites are completely full and completely empty,
 respectively), so from here on we consider initial conditions
\be\label{eq:initcon}
x(0)\in [0,1]^n \setminus \{0_n,1_n\} .
\ee

An important property of the~RFMR
 is
that if~$x_{i+1}(t ) \approx 1$. i.e. site~$i+1$ is ``quite full'' then 
\[
\dot x_{i}\approx \lambda_{i-1}x_{i-1}(1-x_i) \geq 0, 
\]
so the density at site~$i$ increases. Thus, ``traffic jams''
 may gradually evolve behind a full site. 


A calculation shows that the
 Jacobian~$J(x)$ of the~RFMR satisfies~$J=M-D$, where
$D:=\diag(\lambda_n x_n +\lambda_1 (1-x_2),\lambda_1 x_1+\lambda_2(1-x_3),\dots,
\lambda_{n-1} x_{n-1}+\lambda_{n}(1-x_1)  )$, and
\be
 M:=\begin{bmatrix}
						0& \lambda_1 x_1& 0&0& \dots &\lambda_n(1-x_1) \\
						\lambda_1 (1-x_2)&0&\lambda_2 x_2&0& \dots&0\\
						0&\lambda_2(1-x_3)&0            &\ddots&\dots&0\\
						&&\ddots\\
	          0&\dots&  0& \lambda_{n-2}(1-x_{n-1})&0&\lambda_{n-1} x_{n-1} \\
	           \lambda_n x_n&\dots&  0& 0& \lambda_{n-1}(1-x_{n}) & 0  \\
\end{bmatrix}.
\ee
It is straightforward to verify that if for some index~$i$
we have~$x_i(t) \equiv 0$ [$x_i(t) \equiv 1$] on a time interval 
then~$x(t)\equiv 0_n$ [$x(t)\equiv 1_n$] on a time interval. Since we consider initial conditions
as described in~\eqref{eq:initcon}, 
we conclude that~$M(x(t))$ (and thus~$J(x(t))$) 
is irreducible, except perhaps at isolated  points of time. 
Now Theorem~\ref{thm:nsct} implies that the variational
 system~$\dot z(t)=J(x(t))z(t)$ is~{\model}.

\section{Conclusions}\label{sec:conc}

Several interesting studies (see e.g.~\cite{Fusco1990,smith_sign_changes,poin_cyclic}) analyzed 
certain types of
 nonlinear time-varying systems in the form~$\dot y=f (t,y)$ by showing
    that the number of cyclic sign variations in the vector solution~$z(t)$  of  the 
		variational system~$\dot z(t) =J(t,y(t))z(t)$ can only decrease with time. 
This was proved by direct  analysis of  the variational system.

Here, we developed the theoretical framework of such systems 
by
introducing  the notion of a {\model}
 and analyzing  its properties. A~{\model} is a linear time-varying system
whose transition matrix satisfies the strong cyclic variation diminishing property.
We proved that the number
of cyclic sign variations in the vector solution of such a system
can only  decrease with time. We also derived  a necessary and sufficient condition for a system to be~{\model}.
 This generalizes the notion of TPDS analyzed in~\cite{schwarz1970,fulltppaper}.

Our results suggest several 
interesting directions for further research. First, we believe that the 
systematic analysis of~{\model}s presented here can assist
in developing a better understanding  
of   what {\model} of the variational system implies
for the original nonlinear system. 
When the variational system~$\dot z=J z$ is~TPDS then there exists a time~$T\geq 0$ such that~$z(t)\in V$ for all~$t\geq T$
(where~$V$ is the set defined in~\eqref{eq:defvsret}). This means in particular that either~$z_1(t)>0$  for all~$t\geq T$
or~$z_1(t)<0$ for all~$t\geq T$. Since~$z_1=\dot x_1$, we conclude that either~$x_1(t)$ is unbounded or it converges to a limit.
Building on this, Smillie~\cite{smillie} 
 and Smith~\cite{periodic_tridi_smith}  proved powerful stability results for time-invariant and time-varying and periodic
nonlinear systems whose variational  system is~TPDS. However,  Example~\ref{exa:unbou} shows that when the 
variational system  is~{\model}   every~$z_i(t)$ may have an \emph{unbounded} number of sign variations, and thus the
approach used in the~TPDS case cannot be applied. 
 
Nevertheless, the fact that~$s^-_c$ (and~$s^+_c$) is a  discrete-valued Lyapunov function certainly constraints the dynamics.
For example, for a system with a compact state-space it  implies  that 
the  state-space admits a Morse decomposition (see e.g.~\cite{cfsystems}). 
Also, there are certainly special cases where more can be said. 
For example, Remark~\ref{rem:QnotM} shows that if~$s^+_c(z(t))=s^-_c(z(t))=0$
for all~$t\geq T$ then necessarily~$s^+(z(t))=s^-(z(t))=0$, so~$z(t) \in V$ for all~$t\geq T$
 and then   we can apply the same ideas as in the~TPDS case. 
These  issues are  currently under study.

  Other possible topics for further research include   
analyzing  the spectral properties of matrices in~$\Q^+$ (see e.g.~\cite{Fusco1990} for some results on this topic),
 and  studying the properties of the nonlinear~$T$-periodic 
system~$\dot y=f(t,y)$ satisfying~$J(t,y)\in\Q^+$ for all~$t$.

\section*{Acknowledgments}
MM is grateful to  Juergen Garloff  for bringing to his attention the work in Chapter~5 of~\cite{karlin_tp} and to Rola Alsaidi
for helpful comments. 

 
\section*{Biographies}

\noindent {\bf Tsuff Ben Avraham} is an undergraduate student in the EE-Physics program at Tel Aviv University, Israel.

\noindent {\bf Guy Sharon} is an undergraduate student in the EE-Physics program at Tel Aviv University, Israel.

\noindent {\bf Yoram  Zarai}  received the BSc (cum laude), MSc and~PhD degrees in EE
from Tel Aviv University, in 1992, 1998 and~2016 respectively. He is currently a postdoctoral fellow at the Faculty of Engineering, Tel Aviv University. His research interests include modeling and analysis of biological phenomena, systems and control theory and machine learning.    

\noindent {\bf Michael Margaliot}  received the BSc (cum laude) and~MSc degrees in EE from the
Technion-Israel Institute of Technology-in 1992 and 1995, respectively, and the PhD degree
(summa cum laude) from Tel Aviv University in 1999. He was a post-doctoral fellow in the Dept.
of Theoretical Math. at the Weizmann Institute of Science. In~2000, he joined the Department
of EE-Systems, Tel Aviv University, where he is currently a Professor and
Chair. His research interests include the stability analysis of differential inclusions and switched
systems, optimal control theory, fuzzy control, computation with words, Boolean control networks, contraction theory, cooperative systems, 
and systems biology. He is co-author of \emph{New Approaches to Fuzzy Modeling and Control: Design and Analysis}, World
Scientific,~2000 and of \emph{Knowledge-Based Neurocomputing}, Springer,~2009. 
He   served as an Associate Editor
for the journal \emph{IEEE Transactions on Automatic Control} during~$2015-2017$.



 \end{document}